\RequirePackage{ifpdf}
\ifpdf 
\documentclass[pdftex]{sigma}
\else
\documentclass{sigma}
\fi

\usepackage[all]{xy}

\newcommand{\End} {\operatorname{End}\nolimits}
\newcommand{\Hom} {\operatorname{Hom}\nolimits}

\newcommand{\rk}{\operatorname{rk}\nolimits}

\newcommand\CC{{\mathbb C}}

\newcommand\QQ{{\mathbb Q}}
\newcommand\RR{{\mathbb R}}
\newcommand\PP{{\mathbb P}}
\newcommand\ZZ{{\mathbb Z}}
\newcommand\GG{{\mathbb G}}
\newcommand\HH{{\mathbb H}}
\newcommand\NN{{\mathbb N}}

\newcommand{\AAA}{{\mathcal A}}
\newcommand{\CCC}{{\mathcal C}}
\newcommand{\EEE}{{\mathcal E}}
\newcommand{\FFF}{{\mathcal F}}
\newcommand{\HHH}{{\mathcal H}}

\newcommand{\LLL}{{\mathcal L}}
\newcommand{\MMM}{{\mathcal M}}
\newcommand{\NNN}{{\mathcal N}}
\newcommand{\OOO}{{\mathcal O}}
\renewcommand\phi{\varphi}

\newcommand\alp{\alpha}

\newcommand\phiti{\tilde{\phi}}

\newcommand{\eee}{{\boldsymbol e}}

\newcommand{\sss}{{\boldsymbol s}}

\newcommand{\ud}{\mathrm{d}}

\newcommand{\pr}{\operatorname{pr}\nolimits}
\newcommand{\Res}{\operatorname{Res}\nolimits}

\newcommand{\DGal}{\operatorname{DGal}\nolimits}

\newcommand\lra{{\longrightarrow}}

\newcommand\ra{{\rightarrow}}
\newcommand\rar{{\rightarrow}}
\newcommand\equi{{\Longleftrightarrow}}

\renewcommand\bar[1]{\overline{#1}}

\begin{document}

\allowdisplaybreaks

\renewcommand{\PaperNumber}{082}

\FirstPageHeading

\ShortArticleName{Monodromy of a Class of Logarithmic Connections
on an Elliptic Curve}

\ArticleName{Monodromy of a Class of Logarithmic Connections\\ on
an Elliptic Curve}

\Author{Francois-Xavier MACHU} \AuthorNameForHeading{F.-X. Machu}

\Address{Math\'ematiques - b\^{a}t. M2, Universit\'e Lille 1,
F-59655 Villeneuve d'Ascq Cedex, France}
\Email{\href{mailto:xavier.machu@math.univ-lille1.fr}{xavier.machu@math.univ-lille1.fr}}

\ArticleDates{Received March 22, 2007, in f\/inal form August 06,
2007; Published online August 16, 2007}

\Abstract{The logarithmic connections studied in the paper are
direct images of regular connections on line bundles over
genus-$2$ double covers of the elliptic curve. We give an explicit
parametrization of all such connections, determine their
monodromy, dif\/ferential Galois group and the underlying rank-$2$
vector bundle. The latter is described in terms of elementary
transforms. The question of its (semi)-stability is addressed.}


\Classification{14D21; 14H52; 14H60; 32S40}

\section{Introduction}
The Riemann--Hilbert correspondence relates the integrable
logarithmic (or Fuchsian) connections over an algebraic variety
$X$ to the representations of the fundamental group
$\pi_1(X\setminus\{D\})$, where $D$ denotes the divisor of poles
of a connection. Deligne \cite{De} proved its bijectivity, on
condition that $D$ is a f\/ixed normal crossing divisor and the
data on both sides are taken modulo appropriate equivalence
relations. Nevertheless, Deligne's solution is not ef\/fective in
the sense that it does not imply any formulas to compute the
Riemann--Hilbert correspondence. Therefore, it is important to
have on hand a stock of examples that can be solved explicitly.

The authors of \cite{Kor-2,EG} constructed logarithmic connections
of rank $n$ over $\PP^1$ with quasi-permutation monodromy in terms
of theta functions on a ramif\/ied cover of $\PP^1$ of degree~$n$.
Korotkin in~\cite{Kor-1} considers a class of generalized
connections, called connections with constant twists, and
constructs such twisted connections of rank $2$ with logarithmic
singularities on an elliptic curve $E$ via theta functions on a
double cover $C$ of $E$.

In the present paper, we obtain genuine (non-twisted) rank-$2$
connections on $E$ from its double cover $C$ by a dif\/ferent
method, similar to the method applied in \cite{LVU} to the double
covers of $\PP^1$. We consider a genus-$2$ cover $f:C\ra E$ of
degree $2$ with two branch points $p_+$, $p_-$ and a regular
connection $\nabla_\LLL $ on a line bundle $\LLL$ over $C$. Then
the sheaf-theoretic direct image $\EEE=f_* (\LLL)$ is a rank-$2$
vector bundle carrying the connection $\nabla_\EEE :=f_*
(\nabla_\LLL )$ with logarithmic poles at $p_+$ and $p_-$. We
explicitly parameterize all such connections and their monodromy
representations $\rho:\pi_1(E\setminus\{p_{-}, p_{+}\})\ra GL(2,
\CC)$. We also investigate the abstract group-theoretic structure
of the obtained monodromy groups as
 well as their Zariski closures in $GL(2, \CC)$, which are the dif\/ferential Galois groups of the connections $\nabla_\EEE $.

Establishing a bridge between the analytic and algebro-geometric
counterparts of the problem is one of the main objectives of the
paper. We show that the underlying vector bundle $\EEE$ of
$\nabla_\EEE $ is stable of degree $-1$ for generic values of
parameters and identify the special cases where it is unstable and
is the direct sum of two line bundles.

We also illustrate the following Bolibruch--Esnault--Viehweg
Theorem \cite{EV-2}: any irreducible logarithmic connection over a
curve can be converted by a sequence of Gabber's transforms into a
logarithmic connection with same singularities and same monodromy
on a semistable vector bundle of degree $0$. Bolibruch has
established this result in the genus-0 case, in which ``semistable
of degree 0'' means just ``trivial'' \cite{AB}.

We explicitly indicate a Gabber's transform of the above direct
image connection $(\EEE, \nabla_\EEE )$ which satisf\/ies the
conclusion of the Bolibruch--Esnault--Viehweg Theorem. The
importance of results of this type is that they allow us to
consider maps from the moduli space of connections to the moduli
spaces of vector bundles, for only semistable bundles have a
consistent moduli theory. Another useful feature of the elementary
transforms is that they permit to change arbitrarily the degree,
and this enriches our knowledge of the moduli space of connections
providing maps to moduli spaces of vector bundles of dif\/ferent
degrees, which may be quite dif\/ferent and even have dif\/ferent
dimensions (see Remark~\ref{0and-1}).

\looseness=1 All the relevant algebro-geometric tools are
introduced in a way accessible to a non-specialist. One of them is
the usage of ruled surfaces in f\/inding line subbundles of rank-2
vector bundles. This is classical, see~\cite{LN} and references
therein. Another one is the reconstruction of a vector bundle from
the singularities of a given connection on it. Though it is known
as a theoretical method \cite{EV-1,EV-2}, it has not been used for
a practical calculation of vector bundles underlying a given
meromorphic connection over a Riemann surface dif\/ferent from the
sphere. For the Riemann sphere, any vector bundle is the direct
sum of the line bundles $\OOO(k_i)$, and Bolibruch developed the
method of valuations (see~\cite{AB}) serving to calculate the
integers $k_i$ for the underlying vector bundles of connections.
He exploited extensively this method, in particular in his
construction of counter-examples to the Riemann--Hilbert problem
for reducible representations.

Genus-2 double covers of elliptic curves is a classical subject,
originating in the work of Legendre and Jacobi~\cite{J}. We
provide several descriptions of them, based on a more recent
work~\cite{Di}. We determine the locus of their periods (Corollary
\ref{period_locus}), a result which we could not f\/ind elsewhere
in the literature and which we need for f\/inding the image of the
Riemann--Hilbert correspondence in Proposition \ref{imageRHonC}.

\looseness=1 Now we will brief\/ly survey the contents of the
paper by sections. In Section \ref{g2-covers}, we describe the
genus-$2$ covers of elliptic curves of degree $2$ and determine
their periods. In Section \ref{rank-1}, we investigate rank-$1$
connections on $C$ and discuss the dependence of the
Riemann--Hilbert correspondence for these connections on the
parameters of the problem: the period of $C$ and the underlying
line bundle $\LLL$. In Section \ref{Direct_Images}, we compute,
separately for the cases $\LLL=\OOO_C$ and $\LLL\neq\OOO_C$, the
matrix of the direct image connection $\nabla_\EEE $ on $\EEE=f_*
\LLL$. For $\LLL=\OOO_C$, we also provide two dif\/ferent forms
for a scalar ODE of order 2 equivalent to the $2\times2$ matrix
equation $\nabla_\EEE \phi=0$. In Section \ref{monodromy-1}, we
determine the fundamental matrices and the monodromy of
connections $\nabla_\EEE$ and discuss their isomonodromy
deformations. Section \ref{elementary} introduces the elementary
transforms of rank-$2$ vector bundles, relates them to birational
maps between ruled surfaces and states a criterion for
(semi)-stability of a rank-$2$ vector bundle. In Section
\ref{underlying}, we apply the material of Section
\ref{elementary} to describe $\EEE$ as a result of a series of
elementary transforms starting from $\EEE_0=f_* \OOO_C$ and prove
its stability or unstability depending on the value of parameters.
We also describe Gabber's elementary transform which illustrates
the Bolibruch--Esnault--Viehweg Theorem and comment brief\/ly on
the twisted connections of \cite{Kor-1}. In
Section~\ref{monodromy}, we give a description of the structure of
the monodromy and dif\/ferential Galois groups for~$\nabla_\EEE $.

{\bf Terminology.} If not specif\/ied otherwise, a curve will mean
a nonsingular complex projective algebraic curve, which we will
not distinguish from the associated analytic object, a compact
Riemann surface.

\newpage

\section{Genus-2 covers of an elliptic curve}
\label{g2-covers}

In this section, we will describe the degree-2 covers of elliptic
curves which are curves of genus~$2$.

\begin{definition}
Let $\pi:C\rar E$ be a degree-2 map of curves. If $E$ is elliptic,
then we say that $C$ is bielliptic and that $E$ is a degree-2
elliptic subcover of $C$.
\end{definition}

Legendre and Jacobi \cite{J} observed that any genus-2 bielliptic
curve has an equation of the form
\begin{gather}\label{Jacobi}
y^2=c_0x^6+c_1x^4+c_2x^2+c_3\qquad (c_i\in\CC)
\end{gather}
in appropriate af\/f\/ine coordinates $(x,y)$. It immediately
follows that any bielliptic curve $C$ has two elliptic subcovers
$\pi_i:C\rar E_i$,
\begin{gather}
E_1:\ \ y^2=c_0x_1^3+c_1x_1^2+c_2x_1+c_3,\ \ \pi_1:(x,y)\mapsto
(x_1=x^2,y),\ \ \mbox{and}\nonumber\\
E_2:\ \ y_2^2=c_3x_2^3+c_2x_2^2+c_1x_2+c_0,\ \ \pi_2:(x,y)\mapsto
(x_2=1/x^2,y_2=y/x^3).\label{Jacobi-Ei}
\end{gather}
This description of bielliptic curves, though very simple, depends
on an excessive number of parameters. To eliminate unnecessary
parameters, we will represent $E_i$ in the form
\begin{gather}\label{E1-E2}
E_i: y_i^2=x_i(x_i-1)(x_i-t_i) \qquad (t_i\in\CC\setminus \{0,\
1\},\ t_1\neq t_2).
\end{gather}
Remark that any pair of elliptic curves $(E_1,E_2)$ admits such a
representation even if $E_1\simeq E_2$.

We will describe the reconstruction of $C$ starting from
$(E_1,E_2)$ following \cite{Di}. This procedure will allow us to
determine the periods of bielliptic curves $C$ in terms of the
periods of their elliptic subcovers $E_1$, $E_2$.

Let $\phi_i:E_i\rar \PP^1$ be the double cover map
$(x_i,y_i)\mapsto x_i$ ($i=1,2$). Recall that the f\/ibered
product $E_1\times_{\PP^1} E_2$ is the set of pairs $(P_1,P_2)\in
E_1\times E_2$ such that $\phi_1(P_1)= \phi_2(P_2)$. It can be
given by two equations with respect to three af\/f\/ine
coordinates $(x,y_1,y_2)$:
\begin{gather*}
\bar C:=E_1\times_{\PP^1} E_2: \left\{\begin{array}{l}
y_1^2=x(x-1)(x-t_1),\vspace{1mm}\\
y_2^2=x(x-1)(x-t_2).
\end{array}\right.
\end{gather*}
It is easily verif\/ied that $\bar C$ has nodes over the common
branch points $0$, $1$, $\infty$ of $\phi_i$ and is nonsingular
elsewhere. For example, locally at $x=0$, we can choose $y_i$ as a
local parameter on~$E_i$, so that $x$ has a zero of order two on
$E_i$; equivalently, we can write $x=f_i(y_i)y_i^2$ where $f_i$ is
holomorphic and $f_i(0)\neq 0$. Then eliminating $x$, we obtain
that $\bar C$ is given locally  by a single equation
$f_1(y_1)y_1^2=f_2(y_2)y_2^2$. This is the union of two smooth
transversal branches $\sqrt{f_1(y_1)}y_1=\pm \sqrt{f_2(y_2)}y_2$.

Associated to $\bar C$ is its normalization (or desingularization)
$C$ obtained by separating the two branches at each singular
point. Thus $C$ has two points over $x=0$, whilst the only point
of~$\bar C$
 over $x=0$ is the node, which we will denote by the same symbol $0$.
We will also denote by~$0_+$,~$0_-$ the two points of $C$
over~$0$. Any of the functions $y_1$, $y_2$ is a local parameter
at $0_\pm$. In a similar way, we introduce the points $1,\infty\in
\bar C$ and $1_\pm, \infty_\pm\in C$.

\begin{proposition}\label{fp}
Given a genus-$2$ bielliptic curve $C$ with its two elliptic
subcovers $\pi_i:C\rar E_i$, one can choose affine coordinates for
$E_i$ in such a way that $E_i$ are given by the equations
\eqref{E1-E2}, $C$ is the normalization of the nodal curve $\bar
C:=E_1\times_{\PP^1} E_2$, and $\pi_i=\pr_i\circ\nu$, where
$\nu:C\rar\bar C$ denotes the normalization map and $\pr_i$ the
projection onto the $i$-th factor.
\end{proposition}

\begin{proof}
See \cite{Di}.
\end{proof}

It is interesting to know, how the descriptions given by
\eqref{Jacobi} and Proposition \ref{fp} are related to each other.
The answer is given by the following proposition.

\begin{proposition}\label{deg6eq}
Under the assumptions and in the notation of Proposition {\rm
\ref{fp}}, apply the following changes of coordinates in the
equations of the curves $E_i$:
\[
(x_i,y_i)\rar (\tilde x_i, \tilde y_i),\qquad \tilde
x_i=\frac{x_i-t_j}{x_i-t_i},\qquad \tilde y_i=
\frac{y_i}{(x_i-t_i)^2}\sqrt{\frac{(t_j-t_i)^3}{t_i(1-t_i)}},
\]
where $j=3-i$, $i=1,2$, so that $\{i,j\}=\{1,2\}$. Then the
equations of $E_i$ acquire the form
\begin{gather}
E_1:\ \ \tilde y_1^2=\left(\tilde{x_1}-\dfrac{t_2}{t_1}\right)\
\left(\tilde{x_1}-\dfrac{1-t_2}{1-t_1}\right)(\tilde{x_1}-1),\nonumber\\
E_2:\ \ \tilde y_2^2=\left(1-\dfrac{t_2}{t_1}\tilde{x_2}\right)
\left(1-\dfrac{1-t_2}{1-t_1}\tilde{x_2}\right)(1-\tilde{x_2}).\label{eqEi}
\end{gather}
Further, $C$ can be given by the equation
\begin{gather}\label{eqC}
\eta^2=\left(\xi^2-\frac{t_2}{t_1}\right)
\left(\xi^2-\frac{1-t_2}{1-t_1}\right)(\xi^2-1),
\end{gather}
and the maps $\pi_i:C\rar E_i$ by $(\xi,\eta)\mapsto (\tilde x_i,
\tilde y_i)$, where
\[
(\tilde x_1, \tilde y_1)=(\xi^2,\eta),\qquad (\tilde x_2, \tilde
y_2)=(1/\xi^2,\eta/\xi^3).
\]
\end{proposition}

\begin{proof}
We have the following commutative diagram of double cover maps
\begin{gather*}
\xymatrix{ &\ar_{\pi_1}[dl] C \ar_{f}[d] \ar^{\pi_2}[dr]\\
E_1 \ar_{\phi_1}[dr] & \ar_{\tilde{\phi}}[d] {\PP^1} & E_2\ar^{\phi_2}[dl] \\
& {\PP^1} }
\end{gather*}
in which the branch loci of $\phiti$, $\phi_i$, $f$, $\pi_i$ are
respectively $\{t_1,t_2\}$, $\{0,1,t_i,\infty\}$, $\phiti^{-1}
(\{0,1,\infty\})$, $\phi_i^{-1}(t_j)$ ($j=3-i$). Thus the $\PP^1$
in the middle of the diagram can be viewed as the Riemann surface
of the function $\sqrt{\frac{x-t_2}{x-t_1}}$, where $x$ is the
coordinate on the bottom $\PP^1$. We introduce a~coordinate $\xi$
on the middle $\PP^1$ in such a way that $\phiti$ is given by
$\xi\mapsto x$, $\xi^{2}=\frac{x-t_2}{x-t_1}$. Then $C$ is the
double cover of $\PP^1$ branched in the 6 points $\phiti^{-1}
(\{0,1,\infty\})=\big\{\pm1,\pm\sqrt{\frac{1-t_2}{1-t_1}},\pm\sqrt{\frac{t_2}{t_1}}\big\}$,
which implies the equation \eqref{eqC} for $C$. Then we deduce the
equations of $E_i$ in the form \eqref{eqEi} following the recipe
of \eqref{Jacobi-Ei}, and it is an easy exercise to transform them
into \eqref{E1-E2}.
\end{proof}

The locus of bielliptic curves in the moduli space of all the
genus-2 curves is 2-dimensional, hence is a hypersurface.
In~\cite{SV}, an explicit equation of this hypersurface is given
in terms of the Igusa invariants of the genus-2 curves. We will
give a description of the same locus in terms of periods. We start
by recalling necessary def\/initions.

Let $a_1$, $a_2$, $b_1$, $b_2$ be a symplectic basis of
$H_1(C,\ZZ)$ for a genus-2 curve $C$, and $\omega_1$, $\omega_2$ a
basis of the space $\Gamma(C,\Omega^1_C)$ of holomorphic 1-forms
on $C$.

\begin{definition}
Let us introduce the $2\times 2$-matrices $A=(\int_{a_i}\omega_j)$
and $B=(\int_{b_i}\omega_j)$. Their concatenation $\Pi=(A|B)$ is a
$2\times 4$ matrix, called the period matrix of the $1$-forms
$\omega_1$, $\omega_2$ with respect to the basis $a_1$, $a_2$,
$b_1$, $b_2$ of $H_1(C,\ZZ)$. The period of $C$ is the $2\times
2$-matrix $Z=A^{-1}B$. If $A=I$ is the identity matrix, the basis
$\omega_1,\omega_2$ of $\Gamma(C,\Omega^1_C)$ and the
corresponding period matrix $\Pi_0= (I|Z)$ are called normalized.

The period lattice $\Lambda=\Lambda(C)$ is the $\ZZ$-submodule of
rank $4$ in $\Gamma(C,\Omega^1_C)^*$ generated by the $4$~linear
forms $\omega\mapsto \int_{a_i}\omega$, $\omega\mapsto
\int_{b_i}\omega$. A choice of the basis $\omega_i$ identif\/ies
$\Gamma(C,\Omega^1_C)^*$ with $\CC^2$, and $\Lambda$ is then
generated by the $4$ columns of $\Pi$.
\end{definition}

The period $Z_C$ of $C$ is determined modulo the discrete group
Sp$(4,\ZZ)$ acting by symplectic base changes in $H_1(C,\ZZ)$.

{\bf Riemann's bilinear relations.} The period matrix of any
genus-2 curve $C$ satisf\/ies the conditions
\[
Z^t=Z \qquad \mbox{and}\qquad \Im Z >0.
\]

\begin{figure}[t]
\centerline{\includegraphics[width=13cm]{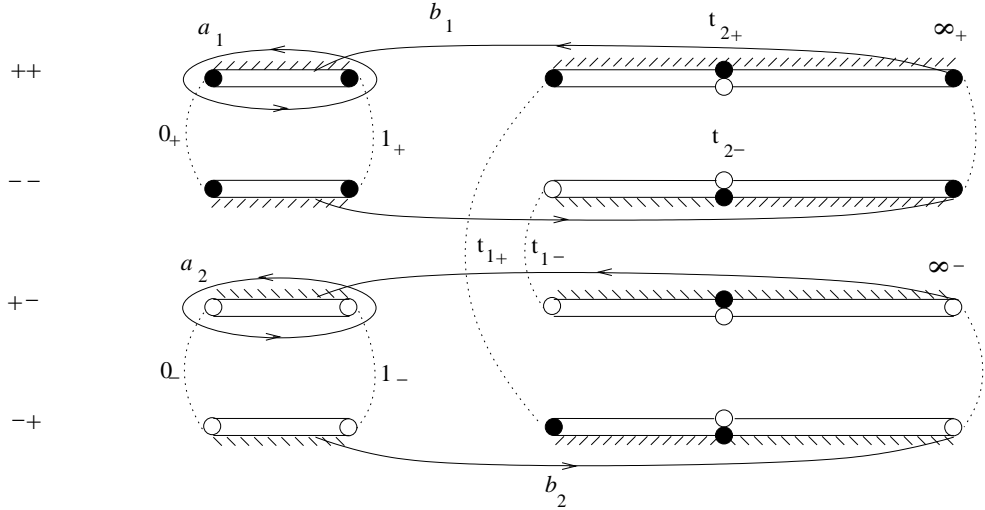}} \caption{The
4 sheets of $C$. The segments of two edges of the cuts are glued
together if they are: (1) situated one under the other, and (2)
hatched by dashes of the same orientation. Thus, the upper edge of
the cut on $\Sigma_{++}$ between $t_1$, $t_2$ is glued to the
lower edge of the cut on $\Sigma_{-+}$ between $t_1$, $t_2$. Four
black points over $t_2$ glue together to give one point $t_{2+}\in
C$, and similarly four white ones give $t_{2-}\in C$. The 4
preimages of each one of the points $0$, $1$, $t_1$, $\infty$ are
glued in pairs, as shown by the colors black/white and by dotted
lines, and give 8 points of $C$ denoted by $0_\pm$, $1_\pm$,
$t_{1\pm}$, $\infty_\pm$.}\label{fig1}
\end{figure}

\looseness=1 To determine the periods of bielliptic curves $C$, it
is easier to use the representation from Proposition \ref{fp}
rather than the standard equation of a genus-2 curve \eqref{eqC}.
This is due to the fact that we can choose $\omega_1=dx/y_1$,
$\omega_2=dx/y_2$ as a basis of the space $\Gamma(C,\Omega^1_C)$
of holomorphic 1-forms on $C$, and the periods of these 1-forms
are easily related to the periods on $E_i$. (Basically,
$(\omega_1, \omega_2)$ can be seen as a basis of eigenvectors of
the action of $(\ZZ/2\ZZ)^{2}$ on~$C$.)

\looseness=1 To f\/ix the ideas, we assume for a while that $t_1$,
$t_2$ are real and $1<t_1<t_2$ (the general case is obtained by a
deformation moving the points $t_i$). $E_i$ can be represented as
the result of gluing two sheets $\Sigma_{i+}$, $\Sigma_{i-}$, which
are Riemann spheres with cuts along the segments $[0, 1]$ and
$[t_i, \infty]$. Then $C$, parameterizing the pairs of points
$(P_1,P_2)$ with $P_i\in E_i$ and with the same $x$-coordinate, is
the result of gluing 4 sheets, which are copies of the Riemann
sphere with cuts along the segments $[0, 1]$ and $[t_1, \infty]$
labelled by $++$, $--$, $+-$, $-+$. For example, the sheet
$\Sigma_{+-}$ is formed by the pairs $(P_1,P_2)$ where $P_1$ lies
on $\Sigma_{1+}$ and $P_2$ on~$\Sigma_{2-}$. Fig.~\ref{fig1} shows
the gluings of the edges of the cuts with the help of hatching and
f\/ixes the choice of the cycles $a_i$, $b_i$. Black points on one
vertical are identif\/ied, the same for the white ones.

\begin{proposition}
Let $C$, $E_1$, $E_2$ be as in Proposition {\rm \ref{fp}}, and
$a_i$, $b_i$ as on Fig.~{\rm \ref{fig1}}. Then the period matrix
of $C$ is
\[
Z_C=\left( \begin{array}{cc}  \frac{1}{2} (\tau_1+\tau_2) &
\frac{1}{2} (\tau_1-\tau_2) \vspace{1mm}\\ \frac{1}{2}
(\tau_1-\tau_2) & \frac{1}{2}(\tau_1+\tau_2)
\end{array}\right),
\]
where $\tau_i$ is the period of $E_i$ with respect to the basis
$\gamma_i=\pi_{i*}(a_1)$, $\delta_i=\pi_{i*}(b_1)$ of
$H_1(E_i,\ZZ)$.
\end{proposition}

\begin{proof} Let $k_i$, $l_i$ be the periods of the dif\/ferential
$dx/y_i$ on $E_i$ along the cycles $\gamma_i$, $\delta_i$
respectively. Take $\omega_i=\pi_i^*(dx/y_i)$ as a basis of
$\Gamma(C,\Omega_C)$. We have
\[
\int_{a_1}\pi_j^*(dx/y_j)=\int_{\pi_{j*}(a_1)}dx/y_j=k_j.
\]
But when calculating the integral over $a_2$, we have to take into
account the fact that a positively oriented loop around a cut on
$\Sigma_{+-}$ projects to a positively oriented loop on
$\Sigma_{2-}$, and the latter def\/ines the cycle $-\gamma_2$ on
$E_{2}$. Thus $\pi_{2*}(a_2)=-\gamma_2$, and the corresponding
period acquires an extra sign:
\[
\int_{a_2}\pi_j^*(dx/y_j)=\int_{\pi_{j*}(a_2)}dx/y_j=(-1)^{j+1}k_j.
\]
The integrals over $b_j$ are transformed in a similar way. We
obtain the period matrix of $C$ in the form
\[
\Pi=\left( \begin{array}{cc|cc} k_1 & k_1 &l_1 & l_1 \\ k_2 & -k_2
&l_2 & -l_2
\end{array}\right).
\]
Multiplying by the inverse of the left $2\times 2$-block and using
the relations $\tau_i=l_i/k_i$, we obtain the result.
\end{proof}

\begin{corollary}\label{period_locus}
The locus $\HHH$ of periods of genus-$2$ curves $C$ with a
degree-2 elliptic subcover is the set of matrices
\[
Z_C=\left( \begin{array}{cc}  \frac{1}{2} (\tau+\tau') &
\frac{1}{2} (\tau-\tau') \vspace{1mm}\\ \frac{1}{2} (\tau-\tau') &
\frac{1}{2}(\tau+\tau')
\end{array}\right)\qquad(\Im\tau>0, \Im\tau'>0)
.\] Equivalently, $\HHH$ is the set of all the matrices of the
form $Z=\left( \begin{array}{cc}  a & b \\ b & a
\end{array}\right)$ ($a,b \in\CC$) such that $\Im Z>0$.
\end{corollary}

\section[Rank-1 connections on $C$ and their monodromy]{Rank-1 connections on $\boldsymbol{C}$ and their monodromy}
\label{rank-1}

We start by recalling the def\/inition of a connection. Let $V$ be
a curve or a complement of a~f\/inite set in a curve $C$. Let
$\EEE$ be a vector bundle of rank $r\geq 1$ on $V$. We denote by
$\OOO_V$, $\Omega^1_V$ the sheaves of holomorphic functions and
1-forms on $V$ respectively. By abuse of notation, we will denote
in the same way vector bundles and the sheaves of their sections.
A {\em connection} on $\EEE$ is a~$\CC$-linear map of sheaves
$\nabla:\EEE\rar\EEE\otimes\Omega^1_V$ which satisf\/ies the
Leibnitz rule: for any open $U\subset V$, $f\in\Gamma(U,\OOO)$ and
$s\in\Gamma(U,\EEE)$, $\nabla(fs)=f\nabla(s)+s\:df$. If $\EEE$ is
trivialized by a basis of sections $\eee=(e_1,\ldots,e_r)$ over
$U$, then we can write $\nabla(e_j)=\sum_ia_{ij}e_i$, and the
matrix $A(\eee)=(a_{ij})$ of holomorphic 1-forms is called the
connection matrix of $\nabla$ with respect to the trivialization
$\eee$. If there is no ambiguity with the choice of a
trivialization, one can write, by abuse of notation, $\nabla
=d+A$.

Given $r$ meromorphic sections $\sss=(s_1,\ldots,s_r)$ which span
$\EEE$ over an open subset, the mat\-rix~$A(\sss)$ def\/ined as
above is a matrix of meromorphic 1-forms on $V$. Its poles in $V$
are called {\em apparent} singularities of the connection with
respect to the meromorphic trivialization $\sss$. The apparent
singularities arise at the points $P\in V$ in which either some of
the $s_i$ are non-regular, or all the $s_i$ are regular but
$s_i(P)$ fail to be linearly independent. They are not
singularities of the connection, but those of the chosen
connection matrix.

In the case when the underlying vector bundle is def\/ined not
only over $V$, but over the whole compact Riemann surface $C$, we
can speak about singularities at the points of $C\setminus V$ of
the connection itself. To this end, choose local trivializations
$\eee_P$ of $\EEE$ at the points $P\in C\setminus V$, and def\/ine
the local connection matrices $A(\eee_P)$ as above, $\nabla
(\eee_P)=\eee_PA(\eee_P)$. The connection~$\nabla$, regular on
$V$, is said to be meromorphic on $C$ if $A(\eee_P)$ has at worst
a pole at $P$ for all $P\in C\setminus V$. If, moreover,
$A(\eee_P)$ can be represented in the form
$A(\eee_P)=B(\tau_P)\frac{d\tau_P}{\tau_P}$, where $\tau_P$ is a
local parameter at $P$ and $B(\tau_P)$ is a matrix of holomorphic
functions in $\tau_P$, then $P$ is said to be a logarithmic
singularity of $\nabla$. A connection is called logarithmic, or
Fuchsian, if it has only logarithmic singularities.

To def\/ine the {\em monodromy} of a connection $\nabla$, we have
to f\/ix a reference point $P_0\in V$ and a basis
$\sss=(s_1,\ldots,s_r)$ of solutions of $\nabla s=0$, $s\in
\Gamma(U,\EEE)$ over a small disc $U$ centered at~$P_0$. The
analytic continuation of the $s_i$ along any loop $\gamma$ based
at $P_0$ provides a new basis
$\sss^\gamma=(s_1^\gamma,\ldots,s_r^\gamma)$, and the monodromy
matrix $M_\gamma$ is def\/ined by $\sss^\gamma =\sss M_\gamma$.
The monodromy matrix depends only on the homotopy class of a loop,
and the monodromy $\rho_\nabla$ of $\nabla$ is the representation
of the fundamental group of $V$ def\/ined by
\[
\rho=\rho_\nabla:\pi_1(V,P_0)\lra GL_r(\CC),\qquad \gamma\mapsto
M_\gamma.
\]

Let now $C=V$ be a genus-2 bielliptic curve with an elliptic
subcover $\phi:C\rar E$. Our objective is the study of rank-2
connections on $E$ which are direct images of rank-1 connections
on $C$. We f\/irst study the rank-1 connections on $C$ and their
monodromy representations.

Let $\LLL$ be a line bundle on $C$ and $e$ a meromorphic section
of $\LLL$ which is not identically zero. Then a connection
$\nabla_\LLL$ on $\LLL$ can be written as $\ud+\omega$, where
$\omega$ is a meromorphic 1-form on $C$ def\/ined by $\nabla_\LLL
(e) =\omega e$. The apparent singularities are simple poles with
integer residues at the points where $e$ fails to be a basis of
$\LLL$. We will start by considering the case when $\LLL$ is the
trivial line bundle $\OOO=\OOO_C$. Then the natural trivialization
of $\LLL$ is $e=1$, and $\omega$ is a regular 1-form. The vector
space $\Gamma(C,\Omega^1_{C})$ of regular 1-forms on $C$ is
2-dimensional; let $\omega_1$, $\omega_2$ be its basis. We can
write $\omega=\lambda_1\omega_1+\lambda_2\omega_2$ with
$\lambda_1$, $\lambda_2$ in $\CC$.

The horizontal sections of $\OOO$ are the solutions of the
equation $\nabla_\OOO \phi=0$. To write down these solutions, we
can represent $C$ as in Proposition \ref{fp} and introduce the
multi-valued functions $z_1=\int \omega_1$ and $z_2=\int
\omega_2$, normalized by $z_1(\infty_+)=z_2(\infty_+)=0$. We
denote by the same symbols $z_1$, $z_2$ the f\/lat coordinates on
the Jacobian $JC=\CC^2/\Lambda$ associated to the basis
$(\omega_1,\omega_2)$ of $\Gamma(C,\Omega^1_{C})$, and $C$ can be
considered as embedded in its Jacobian via the Abel--Jacobi map
$AJ:C\rar JC$, $P\longmapsto\ ((z_{1} (P), z_{2} (P))$ modulo
$\Lambda$.

To determine the monodromy, we will choose $P_0=\infty_+$ and
f\/ix some generators $\alp_i$, $\beta_i$ of $\pi_1(C,\infty_+)$ in
such a way that the natural epimorphism \[ \pi_1(C,\infty_+)\lra
H_1(C,\ZZ)=\pi_1(C,\infty_+)/[\pi_1(C,\infty_+),\pi_1(C,\infty_+)]
\]
 is given by
$\alp_i\mapsto a_i$, $\beta_i\mapsto b_i$.

The following lemma is obvious:

\begin{lemma}\label{obvious}
The general solution of $\nabla_\OOO \phi=0$ is given by $\phi
=ce^{-\lambda_1 z_1-\lambda_2 z_2} $, where $c$ is a~complex
constant. The monodromy matrices of $\nabla_\OOO$ are
$M_{\alp_i}=\exp(-\oint_{a_i} \omega)$,
$M_{\beta_i}=\exp(-\oint_{b_i} \omega)$ ($i=1,2$).
\end{lemma}

Now we turn to the problem of Riemann--Hilbert type: determine the
locus of the representations of $G$ which are monodromies of
connections $\nabla_\LLL$. Since any rank-1 representation $\rho$
of~$G$ is determined by 4 complex numbers
$\rho(\alp_i)$, $\rho(\beta_i)$, we can take $(\CC^*)^4$ for the
moduli space of representations of $G$ in which lives the image of
the Riemann--Hilbert correspondence.

Before solving this problem on $C$, we will do a similar thing on
an elliptic curve $E$. The answer will be used as an auxiliary
result for the problem on $C$.

Any rank-$1$ representation $\rho:\pi_{1}(E)\ra\CC^{*}$ is
determined by the images $\rho(a)$, $\rho(b)$ of the generators $a$,
$b$ of the fundamental group of $E$, so that the space of
representations of $\pi_{1}(E)$ can be identif\/ied with
$\CC^{*}\times\CC^{*}$. We will consider several spaces of
rank-$1$ connections. Let $\CCC(E, \LLL )$ be the space of all the
connections $\nabla:\LLL \ra \LLL \otimes\Omega^{1}_{E}$ on a line
bundle $\LLL $ on $E$. It is non empty if only if $\deg \LLL =0$,
and then $\CCC(E, \LLL )\simeq\Gamma(E, \Omega^{1}_{E})\simeq\CC$.
Further, $\CCC(E)$ will denote the moduli space of pairs $(\LLL ,
\nabla)$, that is, $\CCC(E)=\cup_{[\LLL ]\in J(E)} \CCC(E, \LLL
)$. We will also def\/ine the moduli space  $\CCC$ of triples
$(E_{\tau}, \LLL , \nabla)$, $\CCC=\cup_{\Im\tau>0}
\CCC(E_{\tau})$, and $\CCC_{\rm triv}=\cup_{\Im\tau>0}
\CCC(E_{\tau}, \OOO_{E_{\tau}})$, where $E_{\tau}=\CC/(\ZZ+\ZZ
\tau)$. For any of these moduli spaces, we can consider the
Riemann--Hilbert correspondence map
\[
RH:(E_{\tau}, \LLL , \nabla)\longmapsto (\rho_{\nabla}(a),
\rho_{\nabla} (b)),
\]
where $\rho_{\nabla}$ is the monodromy representation of $\nabla$,
and $(a, b)$ is a basis of $\pi_{1}(E)$ corresponding to the basis
$(1, \tau)$ of the period lattice $\ZZ+\ZZ \tau$. Remark that
$RH\mid_{\CCC(E, \LLL )}$ cannot be surjective by dimensional
reasons. The next proposition shows that  $RH\mid_{\CCC_{\rm
triv}}$ is dominant, though  non-surjective, and that
$RH\mid_{\CCC}$ is surjective.

\begin{proposition}\label{imageRHonE}
In the above notation,
\[
RH(\CCC_{\rm triv})=(\CC^{*}\times\CC^{*}\setminus\{S^{1}\times
S^{1}\}) \cup  \{(1,1)\},\qquad RH(\CCC)=\CC^{*}\times\CC^{*}.
\]
\end{proposition}

\begin{proof}
Let $\nabla=\ud+\omega$ be a connection on an elliptic curve $E$,
where $\omega\in\Gamma(E_{\tau}, \Omega^{1}_{E_{\tau}})$,
$A=\oint_{a} \omega$, $B=\oint_{b} \omega=\tau A$. By analytic
continuation of solutions of the equation $\nabla\phi=0$  along
the cycles in $E$, we obtain $\rho(a)=e^{-A}$ and
$\rho(b)=e^{-\tau A}$. The pair $(-A, -B)=(-A,-\tau A)$ is an
element of ${(0,0)}\cup\CC^{*}\times\CC^{*}$. By setting $z=-A$,
we deduce $RH(\CCC_{\rm triv})=\{(e^{z},
e^{z\tau})\mid(z,\tau)\in\CC\times\HH\}$. The map
$\exp:\CC^{*}\lra\CC^{*}$ is surjective, so for all
$w_1\in\CC^{*}$, we can solve the equation $e^{z}=w_1$, and once
we have f\/ixed $z$, it is possible to solve $e^{\tau z}=w_2$ with
respect to $\tau$ if and only if $(w_1, w_2)\notin S^{1}\times
S^{1} \setminus\{(1,1)\}$. This ends the proof for $RH(\CCC_{\rm
triv})$. The proof for $RH(\CCC)$ is similar to the genus-$2$
case, see Proposition \ref{imageRHonC} below.
\end{proof}

From now on, we turn to the genus-$2$ case. We def\/ine the moduli
spaces $\CCC_2 (C, \LLL )$, $\CCC_2 (C)$, $\CCC_2$, $\CCC{_{2,{\rm
triv}}}$ similarly to the above, so that $\CCC_2(C)=\cup_{[\LLL
]\in J(C)} \CCC_2(C, \LLL )$, $\CCC_2=\cup_{Z\in\HHH}
\CCC_2(C_Z)$, and $\CCC_{\rm triv}=\cup_{Z\in\HHH} \CCC_2(C_Z,
\OOO_{C_Z})$. Here $\HHH$ is the locus of periods introduced in
Corollary \ref{period_locus}, $C_Z$ is the genus-$2$ curve with
period $Z$, $J_2(C)=\CC^{2}/\Lambda$, where $\Lambda\simeq\ZZ^{4}$
is the lattice generated by the column vectors of the full period
matrix $(1\mid Z)$ of $C$. The Riemann--Hilbert correspondence is
the map
\[
RH:(C_Z, \LLL , \nabla)\longmapsto (\rho_{\nabla}(\alp_1),
\rho_{\nabla}(\alp_2) ,\rho_{\nabla}(\beta_1),
\rho_{\nabla}(\beta_2))\in (\CC^{*})^{4},
\]
where the generators $\alp_i$, $\beta_i$ of $\pi_1(C)$ correspond
to the basis of the lattice $\Lambda$.

\begin{proposition}\label{imageRHonC}
In the above notation,
\begin{gather*}
RH(\CCC_{2,{\rm triv}})=\left\{w\in (\CC^{*})^{4}\mid (w_1 w_2, w_3 w_4)\in W, \left(\frac{w_1}{w_2},\frac{w_3}{w_4}\right)\in W \right\},\\
RH(\CCC_2)=(\CC^{*})^{4},
\end{gather*}
where $W$ denotes the locus $RH(\CCC_{\rm triv})$ determined in
Proposition~{\rm \ref{imageRHonE}}.
\end{proposition}

\begin{proof}
Let $\nabla=\ud+\omega$, $\omega\in\Gamma(C_Z, \Omega^{1}_{C_Z}).$
We can consider $C_Z$ in its Abel--Jacobi embedding in $JC$, then
$\omega=\lambda_1\ud z_1+\lambda_2\ud z_2$, where $(z_1, z_2)$ are
the standard f\/lat coordinates on $\CC^{2}/\Lambda$. Therefore,
\[
RH(C_Z,\OOO_Z, \nabla)=\big(e^{\lambda_{1} {z_1}},
e^{\lambda_{2}{z_2}}, e^{\frac{1}{2}(\tau +\tau')\lambda_1
z_1+\frac{1}{2}(\tau-\tau')\lambda_2z_2}, e^{\frac{1}{2}(\tau
-\tau')\lambda_1 z_1+\frac{1}{2}(\tau+\tau')\lambda_2z_2}\big).
\]
Denoting the latter $4$-vector by $w$, we see that $(w_1w_2,
w_3w_4)=(e^{z}, e^{\tau z})$ with $z=\lambda_1z_1+\lambda_2z_2$,
and $(\frac{w_1}{w_2},\frac{w_3}{w_4})=(e^{z'}, e^{\tau'z'})$ with
$z'=\lambda_1z_1-\lambda_2z_2$.

Then Proposition \ref{imageRHonE} implies the answer for
$RH(\CCC_{2,{\rm triv}})$. Now, we will prove the surjectivity of
$RH\mid_{\CCC_2}$. On a genus-2 curve, any line bundle of degree 0
can be represented in the form $\LLL=\OOO_C(P_1+P_2-Q_1-Q_2)$ for
some 4 points $P_i,Q_i\in C$. It is def\/ined by its stalks: for
any $P\in C$, $\LLL_P=\OOO_P$ if $P\not\in \{P_1,P_2,Q_1,Q_2\}$,
$\LLL_{P_i}=\frac{1}{\tau_{P_i}}\OOO_{P_i}$,
$\LLL_{Q_i}=\tau_{Q_i}\OOO_{Q_i}$, where $\tau_P$ denotes a~local
parameter at $P$ for any $P\in C$. This implies that the constant
function $e=1$ considered as a~section of $\LLL$ has simple zeros
at $P_i$ and simple poles at $Q_i$, that is, for its divisor we
can write: $(e)=P_1+P_2-Q_1-Q_2$. According to \cite{A}, any line
bundle of degree 0 admits a connection, and two connections
dif\/fer by a holomorphic 1-form. Hence any connection on $\LLL$
can be written in the form $\nabla=d+\omega$,
$\omega=\nu+\lambda_1dz_1+ \lambda_2dz_2$, where $\nu$ is a
meromorphic 1-form with simple poles at $P_i$, $Q_i$ such that
$\Res_{P_i}\nu=1$, $\Res_{Q_i}\nu=-1$ (these are apparent
singularities of $\nabla$ with respect to the meromorphic
trivialization $e=1$).

We can choose the coef\/f\/icients $\lambda_1$, $\lambda_2$ in such a
way that $\omega$ will have zero $a$-periods. Let us denote the
periods of $\omega$ by $N_i$:
\begin{gather}\label{equaperiod}
N_1=\int_{a_1}\omega,\qquad N_2=\int_{a_2}\omega,\qquad
N_3=\int_{b_1}\omega,\qquad N_4=\int_{b_2}\omega.
\end{gather}
Then $N_1=N_2=0$ by the choice of $\omega$, and
\[
N_{2+j}=2\pi
i\sum_{k}\Res_{s_k}(\omega)\int_{s_0}^{s_k}dz_j,\qquad j=1, 2,
\]
by the Reciprocity Law for dif\/ferentials of 1$^{\rm st}$ and
3$^{\rm rd}$ kinds \cite[Section~2.2]{GH}, where $\sum_{k}s_k$ is
the divisor of poles $(\omega)_\infty$ of $\omega$, and $s_0$ is
any point of $C$. Taking into account that
$(\omega)_\infty=(\nu)_\infty=P_1+P_2+Q_1+Q_2$, $\Res_{P_i}\nu=1$,
$\Res_{Q_i}\nu=-1$, and $z_j(P)=\int_{P_0}^Pdz_j$, we can rewrite:
\[
N_{2+j}=2\pi i [z_j(P_1)-z_j(Q_1)+z_j(P_2)-z_j(Q_2)].
\]
Hence the components of the vector $\frac{1}{2\pi
i}\genfrac{(}{)}{0pt}{0}{N_3}{N_4}$ are the $2$ coordinates on
$JC$ of the class $[\LLL]$ of the line bundle $\LLL$, which is the
same as the divisor class $[P_1+P_2-Q_1-Q_2]$. Now, we can
f\/inish the proof.

Let $(w_i)\in(\CC^*)^4$. Then, we can f\/ind a 1-form $\eta_1$ of
a connection on a degree-0 line bundle~$\LLL _1$ with monodromy
$(1,1,w_3,w_4)$ in choosing $\LLL _1$ with coordinates
$-\frac{1}{2\pi i}(\log w_3, \log w_4)$ on $JC$. In interchanging
the roles of $a$- and $b$-periods, we will f\/ind another 1-form
of connection $\eta_2$ on another degree-0 line bundle $\LLL _2$,
with monodomy $(w_1,w_2,1,1)$. Then $\omega=\eta_1+\eta_2$ is the
form of a connection on $\LLL _1\otimes \LLL _2$ with monodromy
$(w_i)\in(\CC^*)^4$.
\end{proof}

\section{Direct images of rank-1 connections}
\label{Direct_Images}

We will determine the direct image connections
$f_*(\nabla_\LLL)=\nabla_\EEE$ on the rank-$2$ vector bundle
$\EEE=f_{*}\LLL $, where $f:C\rar E$ is an elliptic subcover of
degree 2 of $C$. From now on, we will stick to a representation of
$C$ in the classical form $y^2=F_6(\xi)$, where $F_6$ is a
degree-6 polynomial. We want that $E$ is given the Legendre
equation $ y^2=x(x-1)(x-t), $ but $F_6$ is not so complicated as
in \eqref{eqC}. Of course, this can be done in many dif\/ferent
ways. We will f\/ix for $C$ and $f$ the following choices:
\begin{gather}
f:C=\{y^2=(t'-{\xi^2})(t'-1-{\xi^2})(t'-t-{\xi^2)}\}  \ra
E=\{y^2=x(x-1)(x-t)
)\},\nonumber \\
(\xi,y)   \mapsto (x,y)=(t'-{\xi^2}, y).\label{new_setting}%
\end{gather}

\begin{lemma}
For any bielliptic curve $C$ with an elliptic subcover $f:C\rar E$
of degree $2$, there exist affine coordinates $\xi$, $x$, $y$ on
$C$, $E$ such that $f$, $C$, $E$ are given by \eqref{new_setting} for
some $t,t'\in \CC\setminus \{0, 1\}$, $t\neq t'$.
\end{lemma}

\begin{proof}
By Proposition \ref{fp}, it suf\/f\/ices to verify that the two
elliptic subcovers $E$, $E'$ of the curves~$C$ given by
\eqref{new_setting}, as we vary $t$, $t'$, run over the whole
moduli space of elliptic curves independently from each other.
$E'$ can be determined from \eqref{Jacobi-Ei}. It is a double
cover of $\PP^1$ ramif\/ied at $\frac{1}{t'}$, $\frac{1}{t'-1}$,
$\frac{1}{t'-t}$, $\infty$. This quadruple can be sent by a
homographic transformation to $0$, $1$, $t$, $t'$, hence $E'$ is
given by $y^2=x(x-1)(x-t)(x-t')$. If we f\/ix $t$ and let vary
$t'$, we will obviously obtain all the elliptic curves, which ends
the proof.
\end{proof}

The only branch points of $f$ in $E$ are $p_\pm =(t', \pm y_0)$,
where $y_0=\sqrt{t'(t'-1)(t'-t)}$, and thus the ramif\/ication
points of $f$ in $C$ are $\tilde{p}_{\pm} =(0, \pm y_0)$. In
particular, $f$ is non-ramif\/ied at inf\/inity and the preimage
of $\infty\in E$ is a pair of points $\infty_\pm\in C$.\ \ $E$ is
the quotient of $C$ by the involution $\iota:C\rar C$, called the
Galois involution of the double covering $f$. It is given in
coordinates by $\iota:(\xi, y)\mapsto (-\xi, y)$.

We f\/irst deal with the case when $\LLL $ is the trivial bundle
$\OOO_C$, in which we write $\nabla_\OOO$ instead
of~$\nabla_\LLL$. The direct image $\EEE_0=f_{*}\OOO_C$ is a
vector bundle of rank 2 which splits into the direct sum of the
$\iota$-invariant and anti-invariant subbundles:
$\EEE_0=(f_{*}\OOO_C)^+ \oplus (f_{*}\OOO_C)^-$. The latter
subbundles are def\/ined as sheaves by specifying their sections
over any open subset $U$ of $E$:
\[
\Gamma(U,(f_{*}\OOO_C)^\pm )= \{s\in \Gamma(f^{-1}(U),\OOO_C)\ | \
\iota^*(s)=\pm s\}.
\]
Obviously, the $\iota$-invariant sections are just functions on
$E$, so the f\/irst direct summand $(f_{*}\OOO_C)^+$ is the
trivial bundle $\OOO_E$. The second one is generated over the
af\/f\/ine set $E\setminus \{\infty\}$ by a single generator
$\xi$, one of the two coordinates on $C$. Thus, we can use $(1,
\xi)$ as a basis trivializing $\EEE_0$ over $E\setminus\{\infty\}$
and compute $\nabla=f_*(\nabla_\OOO)$ in this  basis. We use, of
course, the constant function~1 to trivialize $\OOO_C$ and write
$\nabla_\OOO$ in the form
\begin{gather}\label{nablaO}
\nabla_\OOO=d+\omega,\qquad \omega=\nabla_\OOO
(1)=\lambda_1\frac{\ud\xi}{y}+\lambda_2\frac{\xi\ud\xi}{y}.
\end{gather}
Re-writing $\nabla_\OOO (1)=\omega$ in terms of the coordinate
$x=t'-\xi^{2}$, we get:
\[
\nabla_\OOO  (1)=-\frac{\lambda_1}{2(t'-x)}\frac{\ud x}{y} \xi -
\frac{\lambda_2}{2}\frac{\ud x}{y} 1.
\]
Likewise,
\[
\nabla_\OOO (\xi)=-\frac{\lambda_1}{2}\frac{\ud x}{y}1
-\frac{\lambda_2}{2y}\frac{\ud x}{y} \xi-\frac{\ud x}{2(t'-x)}
\xi.
\]
We obtain the matrix of $\nabla=f_*(\nabla_\OOO )$ in the basis
$(1, \xi)$:
\begin{gather}\label{conn_mat}
A=\left( \begin{array}{cc} -\frac{\lambda_2}{2y} \ud x &
-\frac{\lambda_1}{2y} \ud x \vspace{1mm}\\
-\frac{\lambda_1}{2(t'-x)y}\ud x & -\left(\frac{\lambda_2}{2y} +
\frac{1}{2(t'-x)}\right) \ud x
\end{array}\right).
\end{gather}

This matrix has poles at the branch points $p_\pm$ with residues
\begin{gather}\label{residuesA}
{\Res_{p_+} A}=\left( \begin{array}{cc} 0 & 0 \\
\frac{\lambda_1}{2y_0} &  \frac{1}{2}
\end{array}\right), \qquad
{\Res_{p_-} A}=\left( \begin{array}{cc} 0 & 0 \\
-\frac{\lambda_1}{2y_0} &  \frac{1}{2}
\end{array}\right)
.\end{gather}

As the sum of residues of a meromorphic 1-form on a compact
Riemann surface is zero, we can evaluate the residue at
inf\/inity:
\[
\Res_{p_-} A + \Res_{p_+} A =-\Res_{\infty} (A)=\left(
\begin{array}{cc} 0 & 0 \\   0 &  1
\end{array}\right)
.\] It is nonzero, hence $A$ is not regular at $\infty$ and has
exactly 3 poles on $E$. In fact, the pole at $\infty$ is an
apparent singularity due to the fact that $(1, \xi)$ fails to be a
basis of $f_{*}\OOO_C$ at $\infty$, which follows from the
following proposition:

\begin{proposition}
Let $f:C\rar E$ be the bielliptic cover \eqref{new_setting}, and
$\nabla_\OOO=d+\omega$ a regular connection on the trivial bundle
$\OOO_C$ with connection form
$\omega=\lambda_1\frac{\ud\xi}{y}+\lambda_2\frac{\xi\ud\xi}{y}$.
Then the direct image $\nabla=f_*(\nabla_\OOO)$ is a logarithmic
connection on a rank-$2$ vector bundle $\EEE_0$ over $E$, whose
only poles are the two branch points $p_\pm$ of $f$. In an
appropriate trivialization of $\EEE_0$ over
$E\setminus\{\infty\}$, $\nabla$ is given by the connection matrix
\eqref{conn_mat}, and the residues at $p_\pm$ are given by
\eqref{residuesA}.
\end{proposition}

\begin{proof}
If $P\in E$ is not a branch point, then we can choose a small disk
$U$ centered at $P$ such that $f^{-1}(U)$ is the disjoint union of
two disks $U_\pm$. Let $e_\pm$ be a nonzero $\nabla_\OOO$-f\/lat
section of $\OOO_C$ over $U_\pm$. Then $(e_+,e_-)$ is a basis of
$\EEE_0$ over $U$ consisting of $\nabla$-f\/lat sections. This
implies the regularity of $\nabla$ over $U$ (the connection matrix
of $\nabla$ in this basis is zero).

We have shown that the only points where the direct image of a
regular connection might have singularities are the branch points
of the covering. In particular, $\infty$ is not a singularity
of~$\nabla$. The fact that the branch points are logarithmic poles
follows from the calculation preceding the statement of the
proposition.
\end{proof}

At this point, it is appropriate to comment on the horizontal
sections of $\nabla$, which are solutions of the matrix ODE
$d\Phi+A\Phi=0$ for the vector
$\Phi=\left(\!\!\begin{array}{c}\Phi_1\\
\Phi_2\end{array}\!\!\right)$. We remark that the matrix ODE is
equivalent to one scalar equation of second order which we have
not encountered in the literature. It is obtained as follows: the
f\/irst line of the matrix equation gives
\[
\Phi_2=\frac{2\lambda_2}{\lambda_1} \Phi'_1
-\frac{\lambda_2}{\lambda_1} \Phi_1,
\]
where $\Phi_1$, $\Phi_2$ denote the components of a single
$2$-vector $\Phi$, and the prime denotes the derivative with
respect to $x$. The second equation gives:
\[
\Phi'_2=\frac{\lambda_1}{2y(t'-x)} \Phi_1 +
\left(\frac{\lambda_2}{2y} + \frac{1}{2(t'-x)}\right)\Phi_2.
\]
By substituting here $\Phi_2$ in terms of $\Phi_1$, we get one
second order equation for $\Phi_1$. By setting
$y^2=P_3(x)=x(x-1)(x-t)$, we have $y'=\frac{P'_3}{2y}$ and the
dif\/ferential equation for $\Phi_1$ takes the form
\[
\Phi_1''+
\left[\frac{P'_{3}(x)}{2P_{3}(x)}-\frac{\lambda_2}{y}+\frac{1}{2(x-t')}\right]\Phi'_1
+
\left[\frac{\lambda^{2}_{1}}{4P_{3}(x)(x-t')}+\frac{\lambda^{2}_{2}}{4P_3
(x)}-\frac{\lambda_2}{4(x-t')y}\right]\Phi_1=0.
\]
We can also write out the second order dif\/ferential equation for
$\Phi_1$ with respect to the f\/lat coordinate $z=\int \frac{\ud
x}{y}$ on $E$. Now, set up the convention that the prime denotes
$\frac{\ud}{\ud z}$. Then, after an appropriate scaling,
$x=\wp(z)+\frac{t+1}{3}$, $y=\frac{\wp'(z)}{2}$. Let $z_0, -z_0$
be the solutions of $\wp(z)=t'-\frac{t+1}{3}$ modulo the lattice
of periods. Then we have the following equation for $\Phi_1$:
\[
\Phi_1''+\left[-\lambda_2
+\frac{\wp'(z)}{2(\wp(z)-\wp(z_0))}\right]\Phi'_1+
\left[\frac{\lambda^{2}_{2}}{4}+\frac{2\lambda^{2}_{1}-\lambda_2\wp'(z)}{8(\wp(z)-\wp(z_0))}\right]\Phi_1
=0.
\]

We now go over to the general case, in which $\LLL$ is any line
bundle of degree $0$ on $C$ endowed with a regular connection
$\nabla_\LLL$. Then $f_* \LLL=\EEE$ is a vector bundle of rank $2$
on $E$ endowed with a~logarithmic connection $\nabla_\EEE=f_{*}
\nabla_\LLL$. We can represent $\LLL$ in the form
$\LLL=\OOO(\tilde{q}_{1}+ \tilde{q}_{2}- \infty_+ -\infty_-)$ with
$\tilde{q}_{1}=(\xi_1,y_1)$ and $\tilde{q}_{2}=(\xi_2, y_2)$ some
points of $C$. Their images on $E$ will be denoted by $q_i$, or
$(x_i,y_i)$ in coordinates. We will use $1\in\Gamma(\OOO_C)$ as a
meromorphic trivialization of $\LLL$ as in the proof of
Proposition \ref{imageRHonC}. In this trivialization, the
connection form of $\nabla_\LLL$ has simple poles at the 4 points
$\tilde{q}_{i},\infty_\pm$ with residues $+1$ at $\tilde{q}_{i}$
and $-1$ at $\infty_\pm$. It is easy to invent one example of such
a form: $\nu= \frac{1}{2}\left(\frac{y+y_{1}}{\xi-\xi_1}
+\frac{y+y_{2}}{\xi-\xi_2}\right)\frac{\ud\xi}{y}$. Hence the
general form of $\nabla_\LLL$ is as follows:
\begin{gather}\label{eq}\nabla_\LLL=\ud+\omega=\ud+\frac{1}{2}\left(\frac{y+y_{1}}{\xi-\xi_1} +\frac{y+y_{2}}{\xi-\xi_2}\right) \frac{\ud\xi}{y} +\lambda_1 \frac{\ud\xi}{y} +\lambda_2 \frac{\xi\ud\xi}{y}.\end{gather}

We compute $\nabla_\LLL(1)$ and $\nabla_\LLL(\xi)$ and express the
result in the coordinates $(x, y)$ of $E$. This brings us to
formulas for the connection $\nabla_\EEE$ on $E$. We obtain:
\[
\nabla_\LLL(1)
=\omega.1=\left[\frac{1}{2}\left(\frac{y+y_{1}}{\xi(\xi-\xi_1)}
+\frac{y+y_{2}}{\xi (\xi-\xi_2)}\right)+\frac{\lambda_1}{\xi} +
\lambda_2\right]\frac{\xi\ud\xi}{y}.
\]

 Splitting $\frac{1}{\xi-\xi_i}$ into the invariant and anti-invariant parts, we get: \begin{gather*}\nabla_\LLL(\xi)=\xi\nabla_\LLL(1)+ \ud\xi\cdot 1=\Bigg[\frac{1}{2}\left(\frac{(y+y_{1})\xi_1}{\xi^{2}-\xi^{2}_1} +\frac{(y+y_{2})\xi_2}{\xi^{2}-\xi^{2}_2}\right) +\lambda_1 \\
 \phantom{\nabla_\LLL(\xi)=}{} +\frac{1}{2}\left(\frac{y+y_{1}}{\xi^{2}-\xi^{2}_1} +\frac{y+y_{2}}{\xi^{2}-\xi^{2}_2}\right)\xi + \lambda_2\xi+\frac{y}{\xi^{2}} \xi\Bigg]\frac{\xi\ud\xi}{y}.
 \end{gather*}
By using the relations $x=t'-\xi^{2}$, $\ud x=-2\xi\ud\xi$, we
determine the connection $\nabla_\EEE=d+A$, where $A$ is the
matrix of $\nabla_\EEE$ in the basis $(1,\xi)$:
\begin{gather}\label{Conn_Mat}
\left(\!\! \begin{array}{cc} -\frac{1}{2}
(\frac{1}{2}(\frac{y+y_{1}}{x_1-x}
+\frac{y+y_{2}}{x_2-x})+\lambda_2)\frac{\ud x}{y} & -\frac{1}{2}
(\frac{1}{2}(\frac{(y+y_{1})\xi_1}{x_1-x}
+\frac{(y+y_{2})\xi_2}{x_2-x})+\lambda_1)\frac{\ud x}{y}
\vspace{1mm}\\  -\frac{1}{2}(\frac{1}{2}
(\frac{(y+y_{1})\xi_1}{(x_1-x)(t'-x)}
+\frac{(y+y_{2})\xi_2}{(x_2-x)(t'-x)})+\frac{\lambda_1}{t'-x})\frac{\ud
x}{y} & -\frac{1}{2}(\frac{1}{2} (\frac{y+y_{1}}{x_1-x}
+\frac{y+y_{2}}{x_2-x})+\lambda_2+ \frac{y}{t'-x}) \frac{\ud x}{y}
\end{array}\!\!\right)\!.\!\!\!\!
\end{gather}

We compute $\Res_{{p}_{\pm}}A$, where $p_{\pm}$ are the only
singularities of $\nabla_\EEE$:
\begin{gather}\label{Residues}
\Res_{{p}_{\pm}}A=\left( \begin{array}{cc} 0 & 0  \\
\frac{1}{4}(\frac{(y_{1}\pm y_0)\xi_1}{x_1-t'} +\frac{(y_{2}\pm
y_0)\xi_2}{x_2-t'}) \pm\frac{\lambda_1}{2y_0}  &  \frac{1}{2}
\end{array}\right).
\end{gather}

\begin{proposition}\label{deltaE}
Let $f:C\rar E$ be the bielliptic cover \eqref{new_setting},
$\LLL=\OOO(\tilde{q}_{1}+ \tilde{q}_{2}- \infty_+ -\infty_-)$ with
$\tilde{q}_{i}=(\xi_i,y_i)\in C$ ($i=1,2$), and
$\nabla_\LLL=d+\omega$ a regular connection on $\LLL$ with
connection form~$\omega$ defined by \eqref{eq}. Assume that
$\xi_i\neq 0$, that is $\tilde{q}_{i}\neq \tilde p_\pm$. Then the
direct image $\nabla_\EEE=f_*(\nabla_\LLL)$ is a logarithmic
connection on a rank-$2$ vector bundle $\EEE$ over $E$ whose only
poles are the two branch points $p_\pm$ of $f$. In the meromorphic
trivialization of $\EEE$ defined by $(1,\xi)$, $\nabla_\EEE$ is
given by the connection matrix \eqref{Conn_Mat}, and the residues
at $p_\pm$ are given by \eqref{Residues}.
\end{proposition}

Remark that the points $q_i=f(\tilde q_i)=(x_i, y_i)$ are apparent
singularities of $\nabla_\EEE $. We write down the residues of $A$
at these points for future use:
\begin{gather}\label{eqnarray.3.}
{\Res_{q_1}A}=\left( \begin{array}{cc} \frac{1}{2} &
\frac{\xi_1}{2} \vspace{1mm}\\   \frac{1}{2\xi_1 } & \frac{1}{2}
\end{array}\right), \qquad {\Res_{q_2}A}=\left( \begin{array}{cc}
\frac{1}{2} &  \frac{\xi_2}{2} \vspace{1mm}\\   \frac{1}{2\xi_2 }
& \frac{1}{2} \end{array}\right) .
\end{gather}

We can also compute  $\Res_{\infty} A$. First homogenize the
equation of $E$ via the change $x=\frac{x_1}{x_0}$, and
$y=\frac{x_2}{x_0}$. The homogeneous equation is
$x_0x_2^2=x_1^3-(1+t)x_1^2 x_0 +tx_1x_0^2$. Then, setting
$v=\frac{x_0}{x_2}$, $u=\frac{x_1}{x_2}$, we obtain the equation
$v=u^3-(1+t)u^2v+tuv^2$ in the neighborhood of $\infty$. Near
$\infty=(0,0)$, we have $v\sim u^3$, $\frac{\ud x}{y}\sim -2\ud
u$. Therefore,  $\Res_{\infty} A$ is:
\begin{gather}\label{eqnarray.4.}
\Res_{\infty}A=\Res_{u=0}A=\left( \begin{array}{cc} -1 &
-\frac{\xi_1+\xi_2}{2}  \vspace{1mm}\\   0 & -2
\end{array}\right).
\end{gather}

\section{Monodromy of direct image connections}\label{monodromy-1}

\begin{figure}[t]
\centerline{\includegraphics{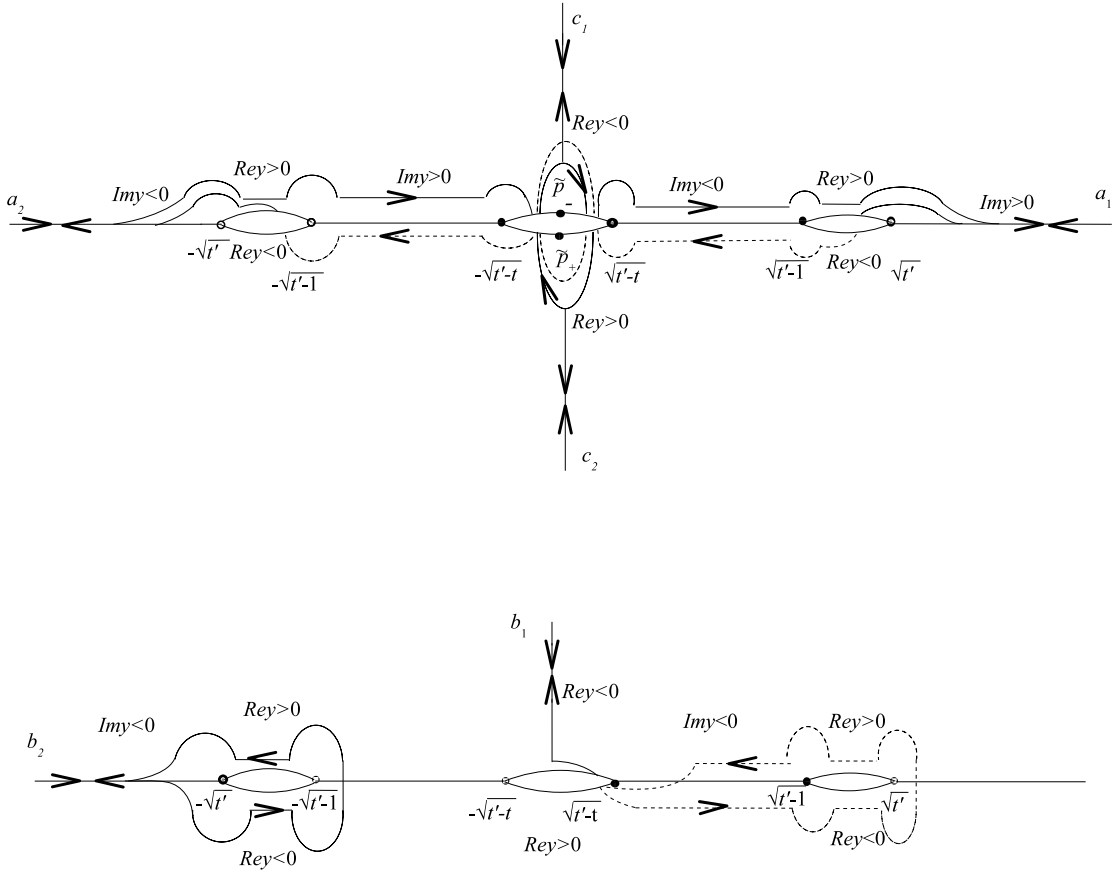}} \caption{Generators of
$\pi_1(C\setminus\{\tilde{p}_\pm\}, \infty_+)$.
 The parts of the arcs represented in solid (resp. dash) lines are on the upper (resp. lower) sheet.}\label{fig3}
\end{figure}

We are using the notation of the previous section. We will
calculate the monodromy of the direct image connections
$\nabla_\EEE$. We will start by choosing generators of the
fundamental group $\pi_1(E\setminus \{p_+,p_-\})$. To express the
monodromy of $\nabla$ in terms of periods of $C$, we will f\/irst
introduce generators $a_i$, $b_i$, $c_i$  for $\pi_1(C\setminus
\{\tilde p_+,\tilde p_-\})$, and then descend some of them to $E$
by applying $f_*$. We choose $\infty_+$ (resp.~$\infty$) as the
reference point on $C$ (resp.~$E$). For this def\/inition, assume
that $t$, $t'$ are real and $1<t<t'$ (for general $t$, $t'$, the
loops $a_i$, $b_i$, $c_i$ are def\/ined up to an isotopy bringing
$t$, $t'$ onto the real axis so that $1<t<t'$). $C$ can be
represented as the result of gluing two copies of the Riemann
sphere along three cuts. We call these copies of the Riemann
sphere upper and lower sheets, and the cuts are realized along the
rectilinear segments $[-\sqrt{t'}$, $-\sqrt{t'-1}]$,
$[-\sqrt{t'-t}, \sqrt{t'-t}]$ and $[\sqrt{t'-1}, \sqrt{t'}]$. The
sheets are glued together in such a way that the upper edge of
each cut on the upper sheet is identif\/ied with the lower edge of
the respective cut on the lower sheet, and vice versa. Let
$\infty_+$ be on the upper sheet, singled out by the condition
$\Im y>0$ when $\xi\in\RR$, $\xi \to +\infty$. This implies that
the values of $\Re{y}$, $\Im{y}$ on $\RR$ are as on Fig.~\ref{fig3},
where the loops $a_i$, $b_i$, $c_i$ generating $\pi_1(C\setminus
\{\tilde p_+,\tilde p_-\})$ are shown. Remark that the loops $c_i$
are chosen in the form $c_i=d_i\tilde c_id_i^{-1}$, where $d_i$ is
a path joining $\infty_+$ with some point close to $\tilde p_\pm$
and $\tilde c_i$ is a small circle around $\tilde p_\pm$ (the
values $i=1,2$ correspond to $\tilde p_+$, $\tilde p_-$
respectively). The paths $d_i$ follow the imaginary axis of the
upper sheet.

Now, we go over to $E$. Set $a=f_*(a_1)$, $b=f_*(b_1)$, and
def\/ine the closed paths running round the branch points $p_\pm$
as follows: $\gamma_i=f(d_i)\tilde\gamma_i f(d_i)^{-1}$, where
$\tilde\gamma_i$ are small circles around~$p_\pm$ running in the
same direction as $f(\tilde c_i)$ (but $f(\tilde c_i)$ makes two
revolutions around $p_\pm$, whilst $\tilde\gamma_i$ only one).

One can verify that the thus def\/ined generators of both\sloppy\
fundamental groups satisfy the relations $[a_1, b_1]c_1[a_2,
b_2]c_2=1$  and $[a,b]\gamma_1\gamma_2=1$ and that the group
morphism
\mbox{$f_*:\pi_1(C\setminus\{\tilde{p}_\pm\},\infty_+)\lra\pi_1(E\setminus\{p_\pm\},\infty)$}
is given by the formulas
\[
f_*(a_1)=a,\quad f_*(b_1)=b,\quad
f_*(a_2)=\gamma^{-1}_{1}a\gamma_1,\quad f_*(b_2)=\gamma^{-1}_{1}
b\gamma_1,\quad f_*(c_i)=\gamma^{2}_{i}\quad (i=1, 2).
\]
As $\nabla_\LLL$ is regular at $\tilde p_\pm$, it has no monodromy
along $c_i$, and this together with the above formulas for $f_*$
immediately implies that the monodromy matrices $M_{\gamma_i}$ of
$\nabla_\EEE$ are of order 2.

We f\/irst assume that $\LLL=\OOO_C$ is trivial, in which case
$\nabla_\LLL$ is denoted $\nabla_\OOO$, and $\nabla_\EEE$ just
$\nabla$. As in the previous section, we trivialize
$\EEE_0=f_{*}(\OOO_C)$ by the basis $(1, \xi)$ over
$E\setminus\{\infty\}$. Splitting the solution $\phi
=e^{-\lambda_1 z_1-\lambda_2 z_2}$ of $\nabla_\OOO\phi=0$ into the
$\iota$-invariant and anti-invariant parts, we represent~$\phi$ by
a $2$-component vector in the basis $(1, \xi)$:
\[
\Phi=\left( \begin{array}{c} e^{-\lambda_2 z_2}\cosh(\lambda_1
z_1)  \vspace{1mm}\\ -\frac{e^{-\lambda_2
z_2}}{\xi}\sinh(\lambda_1 z_1)
\end{array}\right).
\]

We have to complete $\Phi$ to a fundamental matrix
$\mathbf{\Phi}$, and then we can def\/ine the monodromy~$M_\gamma$
along a loop $\gamma$ by $T_\gamma(\mathbf{\Phi})=\mathbf{\Phi}
M_\gamma$, where $T_\gamma$ denotes the analytic continuation
along $\gamma$. We already know the f\/irst column of
$\mathbf{\Phi}$: this is just $\Phi$. Denote it also by
$\mathbf{\Phi_1}$, the column vector
$\left(\begin{array}{c}\Phi_{1,1}\\ \Phi_
{2,1}\end{array}\right)$. It remains to f\/ind $\mathbf{\Phi_2}=
\left(\begin{array}{c}\Phi_{1,2}\\ \Phi_ {2,2}\end{array}\right)$
so that
\[
\mathbf{\Phi}=\left( \begin{array}{cc} \Phi_{1,1} & \Phi_{1,2}
\\  \Phi_{2,1} & \Phi_{2,2}
\end{array}\right)
\] is a fundamental matrix.
By Liouville's theorem, the matrix equation $\mathbf{\Phi'}+
A\mathbf{\Phi}=0$ implies the following scalar equation for
$\Psi=\det\Phi$: $\Psi'+ {\rm Tr}\,(A)\Psi=0$. In our case,
${\rm Tr}\,(A)=-\frac{\lambda_{2}}{y}-\frac{1}{2(t'-x)}$, and we get a
solution in the form: $\Psi=\frac{e^{-2
\lambda_{2}z_{2}}}{\sqrt{t'-x}}=\frac{e^{-2
\lambda_{2}z_{2}}}{\xi}$. Thus we can determine $\mathbf\Phi_2$
from the system:
\begin{gather*}
\cosh(\lambda_1 z_1)\Phi_{2, 2} +\frac{1}{\xi} \sinh(\lambda_1
z_1)\Phi_{1, 2}  =  \frac{e^{-2 \lambda_{2}z_{2}}}{\xi},
\\
\Phi'_{1, 2} =  \frac{1}{2y} (\lambda_2 \Phi_{1, 2} + \lambda_1
\Phi_{2, 2}).
\end{gather*}

Eliminating $\Phi_{2, 2}$, we obtain an inhomogeneous f\/irst
order linear dif\/ferential equation for $\Phi_{1, 2}$. Finally,
we f\/ind:
\[
\mathbf{\Phi_2}=\left( \begin{array}{c} -e^{-\lambda_2
z_2}\sinh(\lambda_1 z_1)  \vspace{1mm}\\ \frac{e^{-\lambda_2
z_2}}{\xi}\cosh(\lambda_1 z_1)
\end{array}\right).
\]

Now we can compute the monodromies of $\nabla_\EEE$ along the
loops $a$, $b$, $\gamma_i$. It is convenient to represent the
result in a form, in which the real and imaginary parts of all the
entries are visible as soon as $t,t'\in \RR$ and $1<t<t'$. Under
this assumption, the entries of the period matrix
$\Pi=((a_{ij}|(b_{ij})$ of $C$ are real or imaginary and can be
expressed in terms of hyperelliptic integrals along the real
segments joining branch points.

Thus reading the cycles of integration from Fig.~\ref{fig3}, we
obtain:
\begin{gather*}
a_{1,1}=-a_{1,2}=2iK, \qquad K=\int_{\sqrt{t'-t}}^{\sqrt{t'-1}} \frac{\ud\xi}{\sqrt{(t'-\xi^2)(t'-1-\xi^2)(t-t'+\xi^2)}}>0,\\
a_{2,1}=a_{2,2}=2iK', \qquad K'=\int_{\sqrt{t'-t}}^{\sqrt{t'-1}} \frac{\xi\ud\xi}{\sqrt{(t'-\xi^2)(t'-1-\xi^2)(t-t'+\xi^2)}}>0,\\
b_{1,1}=-b_{1,2}=-2L, \qquad L=\int_{ \sqrt{t'}}^{\sqrt{t'-1}} \frac{\ud\xi}{\vert{y}\vert}>0,\\
b_{2,1}=b_{2,2}=-2L', \qquad L'=\int_{ \sqrt{t'-1}}^{\sqrt{t'}}
\frac{\xi\ud\xi}{\vert{y}\vert}>0.
\end{gather*}

\begin{proposition}
The monodromy matrices of the connection
$\nabla=f_*(\nabla_\OOO)$, where $\nabla_\OOO$ is the rank-$1$
connection \eqref{nablaO}, are given by
\begin{gather*}
M_{a}=\left( \begin{array}{cc} e^{-2i\lambda_2 K'} \cos(2\lambda_1
K) & -e^{-2i\lambda_2 K'}i\sin(2\lambda_1 K)  \\  -e^{-2i\lambda_2
K'}i\sin(2\lambda_1K) & e^{-2i\lambda_2 K'} \cos(2\lambda_1K)
\end{array}\right),
\\
M_{b}=\left( \begin{array}{cc}
e^{2\lambda_{2}L'}\cosh(2\lambda_1L) &
e^{2\lambda_{2}L'}\sinh(2\lambda_1L)  \\ e^{2\lambda_{2}L'}
\sinh(2\lambda_1L) & e^{2\lambda_{2}L'} \cosh(2\lambda_1L)
\end{array}\right),
\\
M_{\gamma_i}=\left( \begin{array}{cc} 1 & 0 \\ 0 & -1
\end{array}\right) \qquad (i=1,2).
\end{gather*}
\end{proposition}

Now we turn to the general case of nontrivial $\LLL$. Our
computations done in the special case allow us to guess the form
of the fundamental matrix of solutions to $\nabla_\EEE\Phi=0$,
where $\nabla_\EEE=f_*\nabla_\LLL$ and $\nabla_\LLL$ is given by
\eqref{eq} (remark, it would be not so easy to f\/ind it directly
from~\eqref{Conn_Mat}):
\begin{gather*}
\mathbf{\Phi}=\frac{1}{2}\left(\begin{array}{cc}
e^{-\int\omega}+e^{-\int\omega^*} &
e^{-\int\omega}-e^{-\int\omega^*}   \vspace{1mm}\\ \frac{1}{\xi}(
e^{-\int\omega}-e^{-\int\omega^*}) &
\frac{1}{\xi}(e^{-\int\omega}+e^{-\int\omega^*})
\end{array}\right),
\end{gather*}
where $\omega^*:=\iota^*(\omega)$ is obtained from $\omega$ by the
change $\xi\mapsto-\xi$. We deduce the monodromy:

\begin{proposition}
The monodromy matrices of the connection $\nabla_\EEE$ given by
\eqref{Conn_Mat} are the following:
\begin{gather*}
M_{a}=\frac{1}{2}\left( \begin{array}{cc} e^{-N_1}+e^{-N_2} &
e^{-N_1}-e^{-N_2}  \\ e^{-N_1}-e^{-N_2}  & e^{-N_1}+e^{-N_2}
\end{array}\right),
\\
M_{b}=\frac{1}{2}\left( \begin{array}{cc} e^{-N_3}+e^{-N_4} &
e^{-N_3}-e^{-N_4}  \\ e^{-N_3}-e^{-N_4}  & e^{-N_3}+e^{-N_4}
\end{array}\right),
\qquad M_{\gamma_i}=\left( \begin{array}{cc} 1 & 0 \\ 0 & -1
\end{array}\right) \qquad (i=1,2).
\end{gather*}
Here $(N_1, N_2, N_3, N_4)$ are the periods of $\omega$ as defined
in \eqref{equaperiod}.
\end{proposition}

\begin{proof}
By a direct calculation using the observation that the periods of
$\omega^*$ are $(N_2, N_1, N_4$, $N_3)$.
\end{proof}

Here the connection (\ref{Conn_Mat}) depends on $6$ independent
parameters $\tilde{q_1}$, $\tilde{q_2}$, $\lambda_1$, $\lambda_2$,
$t$, $t'$, and the monodromy is determined by the $4$ periods
$N_i$. Hence, it is justif\/ied to speak about the isomonodromic
deformations for this connection. The problem of isomonodromic
deformations is easily solved upon an appropriate change of
parameters. Firstly, change the representation of~$\omega$: write
$\omega=\omega_0+\lambda_1\omega_1+\lambda_2\omega_2$, where
$\omega_0=\nu+\lambda_{10}\omega_1+\lambda_{02}\omega_2$ is chosen
with zero $a$-periods, as in the proof of Proposition
\ref{imageRHonC}, and assume that $(\omega_1, \omega_2)$ is a
normalized basis of dif\/ferentials of f\/irst kind on $C$.
Secondly, replace the $2$ parameters $\tilde{q_1}$, $\tilde{q_2}$
by the coordinates $z_1 [L]$, $z_2 [L]$ of the class of
$L=\OOO(\tilde{q}_{1}+ \tilde{q}_{2}- \infty_+ -\infty_-)$ in
$JC$. Thirdly, replace $(t, t')$ by the period $Z$ of $C$. Then an
isomonodromic variety $N_i=\textrm{const}$ $(i=1,\ldots,4)$ is
def\/ined, in the above parameters, by the equations
\begin{gather*}
\genfrac{(}{)}{0pt}{0}{\lambda_1}{\lambda_2}=\textrm{const},\qquad
\genfrac{(}{)}{0pt}{0}{z_1 [L]}{z_2
[L]}+Z\genfrac{(}{)}{0pt}{0}{\lambda_1}{\lambda_2}=\textrm{const}.
\end{gather*}

Thus the isomonodromy varieties can be considered as surfaces in
the 4-dimensional relative Jacobian $J(\CCC/\HHH)$ of the
universal family of bielliptic curves~$\CCC\rar\HHH$ over the
bielliptic period locus~$\HHH$ introduced in
Corollary~\ref{period_locus}. The f\/iber $C_Z$ of $\CCC$ over a
point $Z\in\HHH$ is a genus-2 curve with period $Z$, and
$J(\CCC/\HHH)\rar\HHH$ is the family of the Jacobians of all the
curves $C_Z$ as $Z$ runs over~$\HHH$. The isomonodromy surfaces
$S_ {\lambda_1,\lambda_2,\mu_1,\mu_2}$ in $J(\CCC/\HHH)$ depend on
4 parameters $\lambda_i$, $\mu_i$. Every isomonodromy surface is a
cross-section of the projection $J(\CCC/\HHH)\rar\HHH$ def\/ined
by
\[
S_{\lambda_1,\lambda_2,\mu_1,\mu_2}= \left\{(Z,[\LLL])\ |\
Z\in\HHH,\ [\LLL]\in JC_Z,\ \genfrac{(}{)}{0pt}{0}{z_1 [L]}{z_2
[L]}=-Z\genfrac{(}{)}{0pt}{0}{\lambda_1}{\lambda_2}
+\genfrac{(}{)}{0pt}{0}{\mu_1}{\mu_2}\right\}.
\]

\section{Elementary transforms of rank-2 vector bundles}
\label{elementary}

In this section, we will recall basic facts on elementary
transforms of vector bundles in the particular case of rank 2, the
only one needed for application to the underlying vector bundles
of the direct image connection in the next section. The impact of
the elementary transforms is twofold. First, they provide a tool
of identif\/ication of vector bundles. If we are given a  vector
bundle $\EEE$ and if we manage to f\/ind a sequence of elementary
transforms which connect $\EEE$ to some ``easy" vector bundle
$\EEE_0$ (like $\OOO\oplus \OOO(-p)$ for a point $p$), we provide
an explicit construction of $\EEE$ and at the same time we
determine, or identify $\EEE$ via this construction. Second, the
elementary transforms permit to change the vector bundle endowed
with a connection without changing the monodromy of the
connection. The importance of such applications is illustrated in
the article~\cite{EV-2}, in which the authors prove that any
irreducible representation of the fundamental group of a Riemann
surface with punctures can be realized by a logarithmic connection
on a {\em semistable} vector bundle of degree 0 (see Theorem
\ref{thm-EV}). On one hand, this is a far-reaching generalization
of Bolibruch's result \cite{AB} which af\/f\/irms the solvability
of the Riemann--Hilbert problem over the Riemann sphere with
punctures, and on the other hand, this theorem gives rise to a map
from the moduli space of connections to the moduli space of vector
bundles, for only the class of {\em semistable} vector bundles has
a consistent moduli theory. We will illustrate this feature of
elementary transforms allowing us to roll between stable,
semistable and unstable bundles in the next section.

Let $E$ be a curve. As before, we identify locally free sheaves on
$E$ with associated vector bundles. Let $\EEE$ be a rank-$2$
vector bundle on $E$, $p$ a point of $E$, $\EEE_{\mid
p}=\EEE\otimes\CC_p$ the f\/iber of $\EEE$ at $p$. Here $\CC_p$ is
the sky-scraper sheaf whose only nonzero stalk is the stalk at
$p$, equal to the 1-dimensional vector space $\CC$. We emphasize
that $\EEE_{\mid p}$ is a $\CC$-vector space of dimension 2, not
to be confused with the stalk $\EEE_p$ of $\EEE$ at $p$, the
latter being a free $\OOO_p$-module of rank 2. Let $e_1$, $e_2$ be
a basis of $\EEE_{\mid p}$. We extend   $e_1$, $e_2$  to sections
of $\EEE$ in a neighborhood of $p$, keeping for them the same
notation. We def\/ine the elementary transforms $\EEE^{+}$ and
$\EEE^{-}$ of $\EEE$ as subsheaves of
$\EEE\otimes\CC(E)\simeq\CC(E)^{2}$, in giving their stalks at all
the points of $E$:
\begin{gather}
\EEE^{-}=elm^{-}_{p,e_2} (\EEE),\qquad \EEE^{-}_{p}=\OOO_p \tau_P e_1+\OOO_p e_2,\nonumber\\
\label{def-elms}\EEE^{+}=elm^{+}_{p,e_1} (\EEE),\qquad \EEE^{+}_{p}=\OOO_p \frac {1}{\tau_P} e_1+\OOO_p e_2,\\
\EEE^{\pm}_{z}=\EEE_z, \qquad \forall\;  z \in
E\setminus\{p\},\nonumber
\end{gather}
where $\tau_P$ denotes a local parameter at $p$. The thus obtained
sheaves are locally free of rank $2$. They f\/it into the exact
triples:
\begin{gather}\label{equa.1.} 0\ra\EEE^{-}\ra\EEE\xymatrix@1{\ar[r]^{\gamma}&} \CC(p)\ra 0,\qquad 0\ra\EEE\ra\EEE^{+}\ra\CC(p)\ra 0.
\end{gather}
Remark, that the surjection $\gamma$ restricted to $\EEE_{\mid p}$
is a projection parallel to the $e_2$ axis; this is the reason for
which we included in the notation of $elm^{-}$ its dependence on
$e_2$. Thus, if we vary~$e_1$, in keeping~$e_2$ (or in keeping the
proportionality class of~$e_2$), the isomorphism class
of~$elm^{-}_{e_2}$ will not change, but it can change if we vary
the proportionality class~$[e_2]$ in the projective line~$\PP
(\EEE_{\mid p})$.

For degrees, we have $\deg \EEE^{\pm}=\deg\EEE\pm 1$. We can give
a more precise version of this equality in terms of the
determinant line bundles: $\det \EEE^{\pm}=\det\EEE (\pm p).$ Here
and further on, given a~line bundle $\LLL$ and a divisor $D=\sum
n_ip_i$ on $E$, we denote by $\LLL(D)$ (``$\LLL$ twisted by $D$'')
the following line bundle, def\/ined as a sheaf by its stalks at
all the points of $E$: $\LLL(D)_z=\LLL_z$ if~$z$~is not among the
$p_i$, and $\LLL(D)_{p_i}= \tau_{p_i}^{-n_i}\LLL_z$. For example,
the regular sections of $\LLL(p)$ can be viewed as meromorphic
sections of $\LLL$ with at most simple pole at $p$, whilst the
regular sections of $\LLL(-p)$ are regular sections of $\LLL$
vanishing at $p$. For the degree of a twist, we have
$\deg\LLL(D)=\deg\LLL+\deg D=\deg\LLL+\sum n_i$, so that $\deg
\LLL(\pm p)=\deg \LLL\pm 1$.

A similar notion of twists applies to higher-rank bundles $\EEE$:
the twist $\EEE(D)$ can be def\/ined either as $\EEE\otimes
\OOO(D)$, or via the stalks in replacing $\LLL$ by $\EEE$ in the
above def\/inition. For degrees, we have $\deg
\EEE(D)=\deg\EEE+\rk \EEE\cdot \deg D$. Coming back to $\rk\EEE=2$
and twisting $\EEE$ by $\pm p$, we obtain some more exact triples:
\[
0\ra\EEE^{+}\ra\EEE(p)\ra\CC(p)\ra 0,\qquad
0\ra\EEE(-p)\ra\EEE^{-}\ra\CC(p)\ra 0.
\]
They are easily def\/ined via stalks, as
$\EEE(p)=\EEE\otimes\OOO_E (p)$ is spanned by $\frac{1}{\tau_P}
e_1$, $\frac{1}{\tau_P} e_2$ at $p$, and $\EEE(-p)$ by
$\tau_Pe_1$, $\tau_Pe_2$.

A basis-free description of elms can be given as follows: Let
$W\subset\EEE_{\mid p}$ be a $1$-dimensional vector subspace. Then
$elm^{-}_{(p, W)}{(\EEE)}$ is def\/ined as the kernel of the
composition of natural maps $\EEE\ra\EEE_{\mid p}\ra\EEE_ {\mid
p}/W$ (here $\EEE_{\mid p}$, $\EEE_ {\mid p}/W$ are considered as
sky-scraper sheaves, i.e.\ vector spaces placed at $p$). The
positive elm is def\/ined via the duality:
\[
elm^{+}_{p, W} (\EEE):=(elm^{-}_{p,W^{\bot}}
(\EEE^{\vee}))^{\vee}.
\]
To set a correspondence with the previous notation, we write:
\[
elm^{+}_{p,e_1} (\EEE)=elm^{+}_{p,\CC e_1} (\EEE),\qquad
elm^{-}_{p,e_2} (\EEE)=elm^{-}_{p,\CC e_2} (\EEE).
\]
One can also def\/ine $elm^{+}$ as an appropriate $elm^{-}$,
applied not to $\EEE$, but to $\EEE(p)$:
\begin{gather}\label{gather.1.} elm^{+}_{p,e_1}=elm^{-}_{p,e_1} (\EEE(p)).\end{gather}

We will now interpret the elementary transforms in terms of ruled
surfaces. For a vector bundle $\EEE$ over $E$, we denote by
$\PP(\EEE)$ the projectivization of $\EEE$, whose f\/iber over
$z\in E$ is the projective line $\PP(\EEE|_z)$ parameterizing
vector lines in $\EEE_{\mid z}$. It has a natural projection
$\PP(\EEE)\rar E$ with f\/ibers isomorphic to $\PP^1$ and is
therefore called a ruled surface. We will see that the elementary
transforms of vector bundles correspond to birational maps between
associated ruled surfaces which split into the composition of one
blowup and one blowdown. The transfer to ruled surfaces allows us
to better understand the structure of $\EEE$, for it replaces all
the line subbundles $\LLL\subset\EEE$ by cross-sections of the
f\/iber bundle $\PP(\EEE)\rar E$; the latter cross-sections being
curves in a surface, we can use the intersection theory on the
surface to study them. As an example, we will give a~criterion of
(semi)stability of $\EEE$ in terms of the intersection theory on
$\PP(\EEE)$.

\begin{figure}[t]
\centerline{\includegraphics[width=13cm]{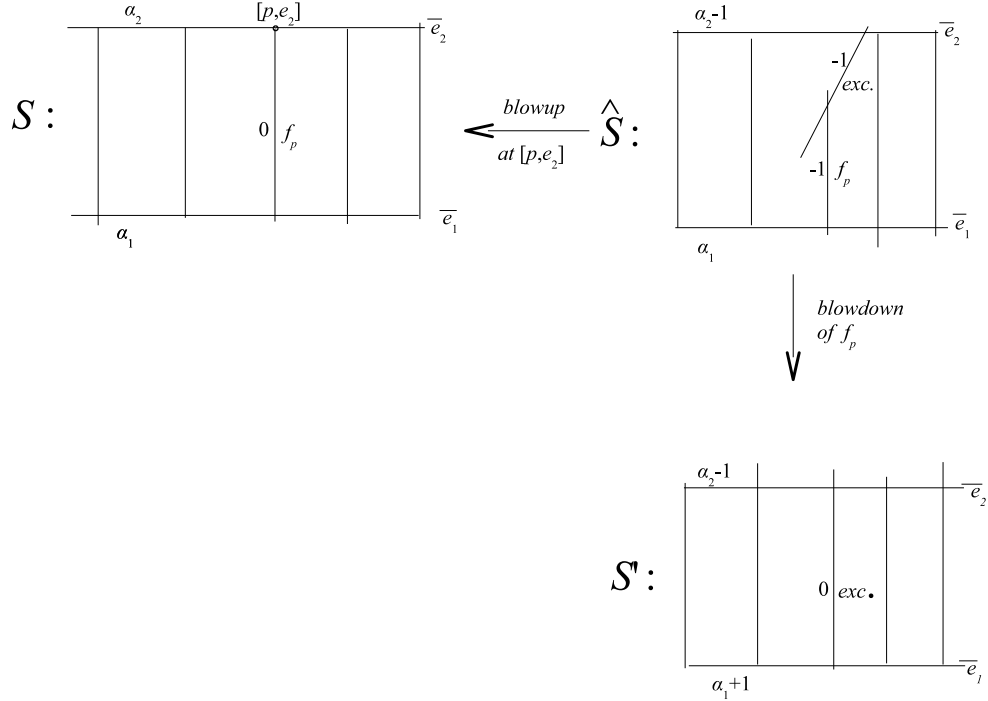}}
\caption{Decomposition of $elm^{-}$ in a blowup followed by a
blowdown.}\label{Fig.6.}
\end{figure}

Let us return to the setting of the description of elms via bases.
We can assume that $e_1$, $e_2$ are rational sections of $\EEE$,
regular and linearly independent at $p$. Let $S=\PP(\EEE)$ and let
$\pi:S\lra E$ be the natural projection. Then $e_1$, $e_2$ def\/ine
two global cross-sections of $\pi$, which will be denoted
$\overline{e_1}$, $\overline{e_2}$. If $\EEE^{-}=elm^{-}_{p, e_2}
(\EEE)$, then the natural map $\EEE^{-}\lra\EEE$ gives rise to the
birationnal isomorphism of ruled surfaces $S\lra
S^{-}=\PP(\EEE^{-})$ which splits into the composition of one
blowup and one blowdown, as shown on Fig.~\ref{Fig.6.}.

Let $f_p$ denote the f\/iber $\pi^{-1}(p)\simeq\PP^1$ of $\pi$; we
keep the same notation for curves and their proper transforms in
birational surfaces. We label some of the curves by their
self-intersection; for example, $
(\overline{e_1}^{2})_S=\alpha_1$, $(f^{2}_p)_S=0$,
$(f^{2}_p)_{\hat{S}}=-1 $. For any vector
$v\in\EEE_{\mid_{p}}\setminus\{0\}$, we denote by $[p, v]$ the
point of $f_p=\PP(\EEE_{\mid_{p}})$ which is the vector line
spanned by~$v$. Remark that the
cross-sections~$\overline{e_1}$,~$\overline{e_2}$ are disjoint in
the neighborhood of $p$ where $e_1$, $e_2$ is a basis of $\EEE$,
but $\overline{e_1}$ can intersect $\overline{e_2}$ at a f\/inite
number of points where $e_1$, $e_2$  fail to generate $\EEE$.

The positive $elm$ has a similar description. Basically, as
$\PP(\EEE)\simeq\PP(\EEE\otimes L)$ for any invertible sheaf $L$
on $E$, we have $\PP(\EEE)\simeq\PP(\EEE(p))$. Hence, in view of
(\ref{gather.1.}), $elm^{+}$ and $elm^{-}$ have the same
representation on the level of ruled surfaces. There exists also
an elegant way to def\/ine $elm^{-}$ in using $\pi:S\lra E$:
\[
elm^{-}_{p, v} (\EEE)=\pi_{*}(I_{S, [p,v]} (1)).
\]
Here, $I_{S, [p,v]}$ is the ideal sheaf of the point $[p, v]$, and
$F(1)$ denotes the twist of a sheaf $F$ by $\OOO_{\PP(\EEE)/E}
(1)$. We have the natural exact triple of an ideal sheaf on
$S=\PP(\EEE)$:
\[
0\ra I_{S, [p,v]}(1)\ra\OOO_{S/E} (1) \ra\CC_{[p, v]}\ra 0.
\]
By a basic property of the tautological sheaf $\OOO_{S/E} (1)$, we
have $\pi_{*}\OOO_{S/E} (1)\simeq \EEE$. By applying~$\pi_{*}$, we
get the exact triple
\[
0\ra \pi_{*}I_{S, [p,v]}(1)\ra\EEE\ra\CC_p\ra 0.
\]
One can prove that in this way we recover the f\/irst exact
triple~(\ref{equa.1.}).

Now, we will say a few words about the (semi)-stability in terms
of ruled surfaces.

\begin{definition}
A rank-$2$ vector bundle on a curve $E$ is stable (resp.\
semistable) if for any line subbundle $\LLL\subset\EEE$,
$\deg\LLL<\frac{1}{2} \deg\EEE$ (resp., $\deg\LLL\leq\frac{1}{2}
\deg\EEE$), or equivalently, if for any surjection $\EEE\rar\MMM$
onto a line bundle $\MMM$, $\deg\MMM>\frac{1}{2} \deg\EEE$ (resp.,
$\deg\MMM\geq\frac{1}{2} \deg\EEE$). A vector bundle is called
unstable if it is not semistable. It is called strictly semistable
if it is semistable, but not stable.
\end{definition}

\begin{definition}
Let $\EEE$ be a rank-$2$ vector bundle on a curve $E$. The index
of the ruled surface  $\pi:S=\PP(\EEE)\lra E$ is the minimal
self-intersection number of a cross-section of $\pi$:
\[
i(S)=\min\{(e)^{2}_{S}\mid e\subset S\textrm{ is a cross-section
of $\pi$}\}.
\]
\end{definition}

The assertion of the following proposition is well-known, see
e.g.~\cite[p.~55]{LN}. For the reader's convenience, we provide a
short proof of it.
\begin{proposition}\label{ss-via-rs}
$\EEE$ is  stable (resp.\ semi-stable) iff $i(S)>0$ (resp.\
$i(S)\geq0$).
\end{proposition}

\begin{proof}
The cross-sections of $\PP(\EEE)\lra E$ are in $1$-to-$1$
correspondence with the exact triples
\begin{gather}\label{equa.2.} 0\ra \LLL _{1}\xymatrix@1{\ar[r]^{\alpha}&}\EEE\xymatrix@1{\ar[r]^{\beta}&} \LLL _{2} \ra 0,\end{gather}
where $\LLL _{1}$, $\LLL _{2}$ are line bundles over $E$. The
cross-section $e$ associated to such a triple is $\PP(\LLL
_1)\subset\PP(\EEE)$. It is the zero locus of
$\pi^{*}\beta\circ\pi^{*}\alpha\in \Hom(\pi^{*}\LLL _1,
\pi^{*}\LLL _2)\simeq H^{0}(S, \pi^{*}(\LLL _2\otimes \LLL
^{-1}_1)).$ Hence, the normal bundle $N_{e/S}$ is isomorphic to
$\LLL _2\otimes \LLL ^{-1}_1$. The stability (resp.\
semi-stability) of $\EEE$ is equivalent to the fact that $\deg
\LLL _1<\deg \LLL _2$ (resp.\ $\deg \LLL _1\leq\deg \LLL _2$) for
any triple (\ref{equa.2.}). As $(e^{2})_S=\deg N_{e/S}=\deg \LLL
_2 - \deg \LLL _1$, this ends the proof.
\end{proof}

We will end this section by two lemmas which help to identify
vector bundles via the geometry of the associated ruled surfaces.

\begin{lemma}\label{direct-sum}
Let $\EEE$ be a rank-$2$ vector bundles over a curve $X$ such that
the associated ruled surface $S=\PP(\EEE)$ has two disjoint
cross-sections $s_1$, $s_2$. Then $\EEE=\LLL_1\oplus\LLL_2$, where
$\LLL_i$ are line subbundles of $\EEE$ corresponding to $s_i$:
$s_i=\PP(\LLL_i)$, $i=1,2$. Further, for the self-intersection
numbers of $s_i$, we have $(s_1)^2=-(s_2)^2=\deg \LLL _2 - \deg
\LLL _1$.
\end{lemma}

\begin{proof}
The f\/irst assertion is obvious, and the second one follows from
the formula for $(e^{2})_S$ in the proof of Proposition
\ref{ss-via-rs}, in taking into account that
$\EEE=\LLL_1\oplus\LLL_2$ f\/its into an exact triple of the form
\eqref{equa.2.}.
\end{proof}

\begin{lemma}\label{elm-on-subb}
Let $\EEE$ be a rank-$2$ vector bundle over a curve $X$, $\LLL$ a
line subbundle of $\EEE$ and $s=\PP(\LLL)$ the associated
cross-section of the ruled surface $S=\PP(\EEE)$. Let $p\in X$,
$[p,v]\in f_p$, where $f_p$ denotes the fiber of $S$ over $p$. Let
$S^\pm=\PP(\EEE^\pm)$, where $\EEE^\pm=elm^\pm_{p,v}$,
$\pi^\pm:S\dasharrow S^\pm$ the natural birational map, $s^\pm$
the proper transform of $s$ in $S^\pm$ under $\pi^\pm$ $($that is,
the closure of $\pi^\pm(s\setminus\{[p,v]\})\,)$, and $\LLL^\pm$
the line subbundle of $\EEE^\pm$ such that $s^\pm=\PP(\LLL^\pm)$.
Then we have:
\begin{enumerate}\itemsep=0pt
\item[(i)] If $[p,v]\in s$, then $(s^\pm)^2_{S^\pm}=(s^2)_S-1 $,
$\deg\LLL^+=\deg\LLL+ 1$, and $\deg\LLL^-=\deg\LLL$. Moreover,
$\LLL^+\simeq\LLL(p)$ and $\LLL^-\simeq\LLL$.

\item[(ii)] If $[p,v]\not\in s$, then
$(s^\pm)^2_{S^\pm}=(s^2)_S+1$, $\deg\LLL^+=\deg\LLL$, and
$\deg\LLL^-=\deg\LLL-1$. Moreover, $\LLL^+\simeq\LLL$ and
$\LLL^-\simeq\LLL(-p)$.
\end{enumerate}
\end{lemma}

\begin{proof}
The formulas for $(s^\pm)^2_{S^\pm}$ follow from the behavior of
the intersection indices as shown on Fig.~\ref{Fig.6.}, and those
for $\deg\LLL^\pm$ are easily deduced directly from the
def\/inition of elementary transforms \eqref{def-elms} by choosing
for $e_1$ or $e_2$ a rational trivialization of $\LLL$.
\end{proof}

\section{Underlying vector bundles of direct image connection}
\label{underlying}

Let us go over again to the setting of Section
\ref{Direct_Images}. Consider f\/irst the case when $\LLL$ is the
trivial bundle, $\LLL=\OOO_C$. The following fact is well known:
\begin{lemma}
Let $f:X\ra Y$ be a finite morphism of smooth varieties of degree
$2$ and $\Delta$ the class of its branch divisor in ${\rm Pic}(Y)$. Then
$\Delta$ is divisible by two in ${\rm Pic}(Y)$, and there exists
$\delta\subset {\rm Pic}(Y)$ such that $2\delta$ is linearly equivalent
to $\Delta$ and $f_{*}\OOO_X=\OOO_Y\oplus\OOO_Y(-\delta)$.
\end{lemma}
\begin{proof}
See \cite[Section 1]{Mu}.
\end{proof}

Applying this lemma to $f:C\ra E$, we f\/ind that
$f_{*}\OOO_C=\OOO_E\oplus\OOO_E(-\delta)$, where $2\delta\simeq
p_+ +p_-$. This property determines $\delta$ only modulo $E[2]$,
but as we saw in Section~\ref{Direct_Images}, $\OOO_E (-\delta)$
is trivialized by a section $\xi$ over $E\setminus\{\infty\}$,
thus $\delta=\infty$ and
$f_{*}\OOO_C=\OOO_E\oplus\OOO_E(-\infty)$. We deduce:

\begin{proposition}
If $\LLL =\OOO_C$, then the direct image connection $\nabla_\EEE
=f_{*}(\nabla_\LLL )$, determined by formula \eqref{Conn_Mat}, is
a logarithmic connection on the vector bundle
$\EEE_0=\OOO_E\oplus\OOO_E(-\infty)$ with two poles at $p_+$,
$p_-$.
\end{proposition}

Let now $\LLL $ be an arbitrary line bundle over $C$ of degree
$0$. By continuity, $\deg f_{*}\LLL =\deg f_{*}\OOO_C=-1$. To
determine $f_{*}\LLL $, we use the following lemma:

\begin{lemma}\label{lemma5.3}\label{oh-la-la}
Let $X$ be a nonsingular curve, $p$ a point in $X$, and $z$ a
local parameter at $p$. Let $\EEE$ be a rank-$2$ vector bundle on
$X$ with a meromorphic connection $\nabla$, regular at $p$. Let
$s_1$, $s_2$ be a pair of meromorphic sections of $\EEE$, linearly
independent over $\CC(X)$ and $\FFF=\langle s_1, s_2\rangle$ the
subsheaf of $\EEE\otimes\CC(X)$ generated by $s_1$, $s_2$ as a
$\OOO_X$-module. Let $\AAA$ be the matrix of $\nabla$ with respect
to the $\CC(X)$-basis $s_1$, $s_2$ and $A={\rm res}_{z=0}\AAA$.
Assume that $\EEE_{\mid p}$ has a basis $v_1$, $v_2$ consisting of
eigenvectors of $A$.  Then $v_1$, $v_2$ extend to a basis of the
stalk $\EEE_p$, the corresponding eigenvalues $n_1$, $n_2$ of $A$
are integers and we have the following relations between the
stalks of subsheaves of $\EEE\otimes\CC(X)$ at $p$:
\begin{gather*}
\EEE_p=\langle v_1, v_2\rangle,\qquad \FFF_p=\langle z^{n_{1}} v_1, z^{n_{2}} v_2\rangle,\\
{\rm if} \ \ n_1=1,\  \ n_2=0,\quad  \EEE_{p}=elm^{+}_{p, v_1} (F_p),\\
{\rm if} \ \ n_1=-1,\ \ n_2=0,\quad \EEE_{p}=elm^{-}_{p, v_2} (F_p),\\
{\rm if} \ \ n_1=n_2,\quad \EEE_p=(F(n_1))_p.
\end{gather*}
\end{lemma}
\begin{proof}
Straightforward.
\end{proof}

Let us apply this lemma to the connection $\nabla_\EEE $, given by
formula (\ref{Conn_Mat}) in  the basis $(1, \xi)$. We have:
$\EEE_{p}\neq \FFF_p\iff p \in\{q_1, q_2, \infty\},$
\[
 A_{i}={\Res_{q_i}\AAA}=\left( \begin{array}{cc} \frac{1}{2} &  \frac{\xi_1}{2} \\   \frac{1}{2\xi_1 } & \frac{1}{2}
\end{array}\right) \quad (i=1, 2),\qquad A_{\infty}=\Res_{\infty}{\AAA}=\left( \begin{array}{cc} -1 & -\frac{\xi_1+\xi_2}{2}  \\   0 & -2
\end{array}\right).
\]
We list the eigenvectors $v^{(i)}_{j}$, $v^{(\infty)}_{j}$
together with the respective eigenvalues for the matrices $A_i$,
$A_{\infty}$:
\begin{gather*}
v^{(i)}_{1}=\left( \begin{array}{c} - \xi_i \\  1
\end{array}\right), \quad \eta^{i}_{1}=0, \qquad \ v^{(i)}_{2}=\left( \begin{array}{c} \xi_i \\  -1
\end{array}\right), \quad \eta^{i}_{2}=-1, \\
 \ v^{(\infty)}_{1}=\left( \begin{array}{c} 1 \\  0
\end{array}\right),\quad \eta^{\infty}_{1}=1, \qquad \ v^{(\infty)}_{2}=\left( \begin{array}{c} -\frac{\xi_1+\xi_2}{2} \\  1
\end{array}\right),\quad \eta^{\infty}_{2}=-2.
\end{gather*}

Applying Lemma \ref{lemma5.3} (twice at $\infty$), we obtain the
following corollary:
\begin{corollary}\label{two-elms}
Let $\LLL=\OOO_C(\tilde q_1+\tilde q_2-\infty_+-\infty_-)$,
$\tilde q_i=(\xi_i, y_i)$, $q_i=f(\tilde q_i)$, $i=1,2$, as in
Proposition \ref{deltaE}, and let $v_j^{(i)}$ be the eigenvectors
of $A_i$ as above. Then
\[
\EEE=elm^{+}_{q_1, v_2^{(1)}}elm^{+}_{q_2, v_2^{(2)}}
(\EEE_0(-\infty)).
\]
\end{corollary}

\begin{remark}
Note that though the sheaf-theoretic direct image $f_*\LLL$ does
not depend on the choice of a connection $\nabla_\LLL$ on $\LLL$,
our method of computation of $f_*\LLL$, given by Corollary
\ref{two-elms}, uses the direct image connection
$\nabla_\EEE=f_*\nabla_\LLL$ for some $\nabla_\LLL$.
\end{remark}

\begin{proposition}\label{generic-case}
For generic $\LLL \in {\rm Pic}(C)$, the rank-$2$ vector bundle
$\EEE$ is stable.
\end{proposition}

\begin{proof}
Starting from the ruled surface
$S_0=\PP(\EEE_0)=\PP(\OOO_E\oplus\OOO_E(-\infty))$, we apply two
elementary transforms $S_0\ra S_1\ra S_2=\PP(\EEE)$, and we have
to prove that any cross-section of $S_2$ has strictly positive
self-intersection, provided that $\tilde q_i=(\xi_i, y_i)$ are
suf\/f\/iciently generic. For a rational section $s$ of $\EEE_0$,
let us denote by $\overline s$ the associated  cross-section of
$S_0$. $S_0$ is characterized by the existence of two
distinguished sections  $\overline s_1$, $\overline s_2 $
associated to $s_1=1$, $s_2=\xi$ with self-intersections
$\overline s_1 ^2=-1$, $\overline s_2 ^2=1$, and we have the
relations $\overline s_1 \overline s_2=0$, $\overline s_2\sim
\overline s_1 + f_{\infty}$, where $f_p=\pi^{-1}(p)$ is the
f\/iber of the structure projection $\pi:S_0\lra E$. When there is
no risk of confusion, we will keep the same notation for curves
and their proper transforms in birational surfaces. Any
cross-section~$\overline s$ is linearly equivalent to $\overline
s_1 + f_{p_1}+\dots +f_{p_r}$ for some points $p_1,\ldots, p_r$ in
$E$, and  $\overline s^2=2r+1$. In particular, $i(S_0)=-1$,
attained on $\overline s_1$. Remark that $\overline s_0$ is rigid,
whilst $\overline s_1$ moves in a pencil $|\overline s_1+
f_{\infty}|$. Let us apply $elm^{+}_{q_1, v_2^{(1)}}$. First, we
blow up $P_1=[q_1, v_2^{(1)}]$. Let $\overline e_1$ be the
corresponding $(-1)$-curve and $\hat{S_0}$ the blown up surface.
For the self-intersection numbers of the cross-sections, we have
the following relations: $(\overline s^2)_{\hat{S_0}}=(\overline
s^2)_{S_0}$ if $P_1\notin\overline s$ and $(\overline
s^2)_{\hat{S_0}}=(\overline s^2)_{S_0} -1$ if $P_1\in \overline
s$. Hence, $\hat{S_0}$ has only one cross-section for each one of
the self-intersection numbers $-1$, $0$, and  $(\overline
s^2)_{S_0} \geq 1$ for all the other cross-sections. The
cross-section with self-intersection $-1$ is $\overline s_1$ and
the one with self-intersection $0$ is the proper transform of the
unique member $\overline s_{P_1}$ of the pencil $\mid \overline
s_1+ f_\infty \mid$ on $S_0$ going through $P_1$, see
Fig.~\ref{Fig.7.}. The next step is the blowdown of
$f_{q_{1}}\subset\hat{S_0}$. The self-intersection number of all
the cross-sections of~$\hat{S_0}\ra E$  that meet $f_{q_{1}}$ goes
up by $1$. We conclude that $S_1=\PP(elm^{+}_{q_1, v_2^{1}}
(\EEE_0))$ has two cross-sections~$\overline s_{P_1}$, $\overline
s_1$ with square~$0$, and $(\overline s^{2})_{S_1}\geq 2$ for any
other cross-section of~$S_1$. In the language of vector bundles,
this means that $\EEE_1$ is the direct sum of two line bundles of
degree $0$. More precisely, $\EEE_1=\OOO_E\oplus\OOO_E (q_1
-\infty)$ by Lemma~\ref{direct-sum}, the f\/irst summand
corresponding to $\overline s_1$ and the second one to $\overline
s_{P_1}$.

\begin{figure}[t]
\centerline{\includegraphics[width=9cm]{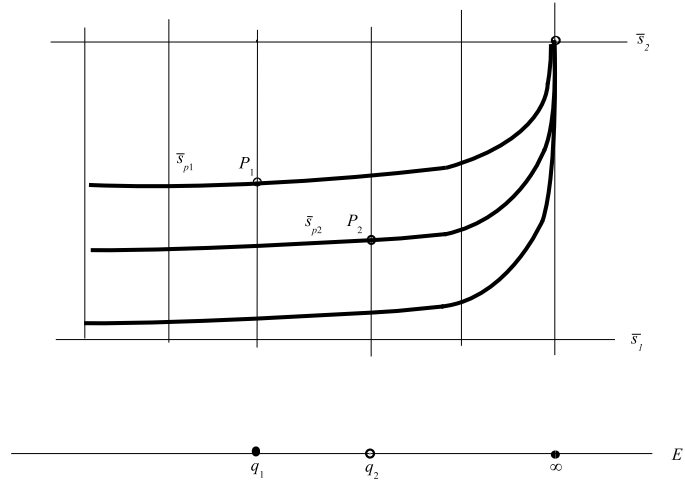}} \caption{The
ruled surface $S_0$. The pencil $\mid \bar s_1+f_\infty \mid$ has
a unique member passing through $P_i$ for each $i=1, 2$.}
\label{Fig.7.}
\end{figure}

The second elementary transform is performed at $P_2\in S_1$. As
$P_2\notin\overline s_{P_1}\cup\overline s_1$, the minimal
self-intersection number of a cross-section in $S_1$ passing
through $P_2$ is $2$. The elementary transform decreases by $1$
the self-intersection of such cross-sections and increases by $1$
the self-intersection of all other cross-sections (Lemma
\ref{elm-on-subb}). Hence, $i(S_2)=1$, the value attained on many
cross-sections, for example, $\overline s_{P_2}$, $\overline
s_{P_1}$, $\overline s_1$. This ends the proof.
\end{proof}

\begin{theorem}[Atiyah, \cite{A}]
For any line bundle $\NNN $ of odd degree over an elliptic curve
$E$, there exists one and only one stable rank-$2$ vector bundle
on $E$ with determinant $\NNN $.
\end{theorem}

Using Atiyah's theorem in our case, we have $\deg \NNN =-1$ , so
that $\NNN $ can be represented in the form $\NNN =\OOO_E (-q)$
for some $q\in E$. $\EEE$ is obtained as the unique non-trivial
extension of vector bundles:
\[
0\ra\OOO_E(-q)\ra\EEE\ra\OOO_E\ra 0.
\]
Moreover, the correspondence $\EEE\leftrightarrow q$ identif\/ies
the moduli space $\MMM^{s}_{E} (2, -1)$ of rank-$2$ stable vector
bundles of degree $-1$ over $E$ with $E$ itself. We deduce:
\begin{corollary}\label{map-to-E}
Under the above identification $\MMM^{s}_{E} (2, -1)\simeq E$, the
rational map:
\begin{gather*}
f:JC  \dashrightarrow  \MMM^{s}_{E} (2, -1),
\\
\LLL =\OOO_C(\tilde{q_1} +\tilde{q_2}-\infty_+ -\infty_-)  \mapsto
f_{*}(\LLL )
\end{gather*}
can be given by
\[
[\tilde{q_1} +\tilde{q_2}-\infty_+ -\infty_-]  \mapsto
[q_1+q_2-2\infty].
\]
\end{corollary}

Now we go over to the nongeneric line bundles $\LLL$. The direct
image $f_{*}\LLL $ can be unstable for special $\LLL $. This may
happen when either the argument of Proposition \ref{generic-case}
does not work anymore, or when formulas
(\ref{Conn_Mat})--(\ref{eqnarray.4.}) are not valid. We list the
cases which need a separate analysis in the next proposition.

\begin{proposition} Let
$\LLL =\OOO_C(\tilde{q_1} +\tilde{q_2}-\infty_+ -\infty_-)$,
$\EEE=f_*(\LLL)$, $\EEE_0=f_*\OOO_C$, as above. Whenever
$\tilde{q_i}$ is finite, it will be represented by its
coordinates: $\tilde{q_i}=(\xi_i,y_i)$. The following assertions
hold:
\begin{enumerate}\itemsep=0pt
\item[(a)] If $\tilde{q_1} +\tilde{q_2}$ is a divisor in the
hyperelliptic linear series $g^{1}_{2}(C)$ (that is $\xi_1=\xi_2$,
$y_1=-y_2$, or
$\{\tilde{q}_1,\tilde{q}_2\}=\{\infty_+,\infty_-\}$), then
$\EEE\simeq\EEE_0$, and hence $\EEE$ is unstable.

\item[(b)] If $\tilde{q_1}=\tilde{q_2}\neq \infty_\pm$, then
$\EEE\simeq \OOO_E{(-\infty)}\oplus\OOO_E({2q_1-2\infty})$ is
unstable.

\item[(c)] If $\tilde{q_i}={\infty_{\pm}}$ for at least one value
$i\in \{1,2\}$, then $\EEE\simeq
\OOO_E(-2\infty+q_{3-i})\oplus\OOO_{E}$ is unstable.

\item[(d)] If $\tilde{q_i}=\tilde{p}_\pm$ for exactly one value
$i\in \{1,2\}$, then $\EEE$ is a stable bundle of  degree $-1$
with $\det\EEE\simeq \OOO_E(q_{3-i}+p_\pm-3\infty)$.
\end{enumerate}
\end{proposition}

\begin{proof}
(a) In this case, $s_{p_1}=s_{p_2}$, $\tilde{q_1}
+\tilde{q_2}\sim{\infty_+}+{\infty_-}$, then $\LLL \simeq\OOO_C$,
and $\EEE\simeq\EEE_0$.

(b) Let, for example, $i=1$. Then $\LLL =\OOO_C(\tilde{q_1}
+\tilde{q_2}-\infty_+ -\infty_-)$ degenerates to $\LLL
=\OOO_C(2\tilde{q_1}-\infty_+ -\infty_-)$. In this case, the
matrix of a regular connection on $\LLL $ in the rational
basis~$1$ of $\LLL =\OOO_C(2\tilde{q_1}-\infty_+
-\infty_-)\hookrightarrow\OOO_C(2\tilde{q_1})\hookleftarrow\OOO_C\ni
1$ is a rational $1$-form with residues~$2$ at~$\tilde{q_1}$ and
$-1$ at points~$\infty_{\pm}$. Such a $1$-form can be written by
the same formula $\omega=\frac{1}{2}\big(\frac{y+y_{1}}{\xi-\xi_1}
+\frac{y+y_{2}}{\xi-\xi_2}\big) \frac{\ud\xi}{y} +\lambda_1
\frac{\ud\xi}{y} +\lambda_2 \frac{\xi\ud\xi}{y}$, as in the
general case, but now we substitute $\xi_2=\xi_1$, $y_2=y_1$ in
it:
\[
\omega=\frac{y+y_{1}}{\xi-\xi_1}\frac{\ud\xi}{y} +\lambda_1
\frac{\ud\xi}{y} +\lambda_2 \frac{\xi\ud\xi}{y}.
\]
Assume that $\xi_1=\xi_2\neq 0$; the case when $\tilde q_1=\tilde
q_2=\tilde p_\pm$ should be treated separately as in~(d). Then we
get the matrix $A$ of $\nabla_\EEE $ in substituting
$\xi_2=\xi_1$, $y_2=y_1$, $x_2=x_1$ into formulas
(\ref{Conn_Mat})--(\ref{eqnarray.4.}). We obtain the following
residues:
\begin{gather*}
\Res_{p\pm}A=\left( \begin{array}{cc} 0 & 0  \\
\frac{1}{2}\frac{(y_1\pm
y_{0})\xi_1}{(x_1-t')}\pm\frac{\lambda_1}{2y_0}  &  \frac{1}{2}
\end{array}\right), \qquad \Res_{q_1}A=\left( \begin{array}{cc} 1 &  \xi_1 \\   \xi_1  & 1
\end{array}\right),\\
\Res_{\infty}A=\Res_{u=0}A=\left( \begin{array}{cc} -1 & -\xi_1
\\   0 & -2
\end{array}\right).
\end{gather*}

As in the proof of Proposition \ref{generic-case}, we can describe
$\EEE$ as the result of two  successive positive elm's applied to
$\EEE_0(-\infty)$. In contrast to the general case, considered in
Lemma \ref{oh-la-la}, the second elm has for  its center the point
$\tilde{P_1}=\overline{s}_{P_1} \cap\tilde{f_{q_1}}\subset S_1$,
where $\tilde{f_{q_1}}$ is the f\/iber of $S_1\ra E$ over~$q_1$.
As $(\overline{s}_{P_1})^2_{S_1}=0$, the resulting surface $S_2$
has a cross-section with self-intersection $-1$, thus $i(S_2)=-1$,
and consequently $\EEE$ is unstable. Applying Lemmas
\ref{direct-sum} and \ref{elm-on-subb}, we can identify it with
$\OOO_E{(-\infty)}\oplus\OOO_E({2q_1-2\infty})$.

(c) Let, for example, $\tilde{q_2}={\infty_{-}}$. Then $\LLL $
degenerates to $\OOO_C(\tilde{q_1}-{\infty_{+}})$, and we can
again write the connection in the same way as in the previous
case. $\EEE$ is obtained from $\EEE_0(-\infty)$ by $2$ positive
elms. From Lemmas \ref{direct-sum} and \ref{elm-on-subb}, we
deduce that $\EEE\simeq\OOO_E(-2\infty+q_1)\oplus\OOO_{E}$.

(d) One of the points $\tilde{p}_{\pm}$ collides with
$\tilde{q_i}$. This corresponds to $\xi_{i}= 0$. So, we assume
that $\xi_{2}=0$, $\xi_{1}\neq 0$ ($\tilde{q_2}=\tilde{p}_+$).
Hence, the $1$-form of the connection can be written as follows:
\[
\omega=\frac{1}{2}\left(\frac{y+y_0}{\xi}+\frac{y+y_{1}}{\xi-\xi_1}\right)\frac{\ud\xi}{y}
+\lambda_1 \frac{\ud\xi}{y} +\lambda_2 \frac{\xi\ud\xi}{y}.
\]
The matrix $A$ is given by
\begin{gather*}
A=\left( \begin{array}{cc} -\frac{1}{2}
\left(\frac{1}{2}\frac{y+y_{1}}{x_1-x}
+\frac{y+y_{0}}{t'-x}+\lambda_2\right)\frac{\ud x}{y} &
-\frac{1}{2}\left(\frac{1}{2}
\frac{(y+y_{1})\xi_1}{x_1-x}+\lambda_1\right)\frac{\ud x}{y}
\vspace{1mm}\\  -\frac{1}{2}\left(\frac{1}{2}
\frac{(y+y_{1})\xi_1}{(x_1-x)(t'-x)}+\frac{\lambda_1}{t'-x}\right)\frac{\ud
x}{y} & -\frac{1}{2}\left(\frac{1}{2}\left(\frac{y+y_{1}}{x_1-x}
+\frac{y+y_{0}}{t'-x}\right)+\lambda_2 +\frac{y}{t'-x}\right)
\frac{\ud x}{y}
\end{array}\right).
\end{gather*}
Its residues at f\/inite points are:
\begin{gather*}
\Res_{p_+}A=\left( \begin{array}{cc} \frac{1}{2} & 0
\vspace{1mm}\\
\frac{1}{4}\frac{(y_0+y_{1})\xi_1}{(x_1-t')}+\frac{\lambda_1}{2y_0}
&  1
\end{array}\right),
\qquad \Res_{p_-}A=\left( \begin{array}{cc} 0 & 0  \vspace{1mm}\\
\frac{1}{4}\frac{(y_1-y_0)\xi_1}{(x_1-t')}-\frac{\lambda_1}{2y_0}
&  \frac{1}{2}
\end{array}\right),
\\
\Res_{q_1}A=\left( \begin{array}{cc} \frac{1}{2} &
\frac{\xi_1}{2} \vspace{1mm}\\   \frac{1}{2\xi_1 } & \frac{1}{2}
\end{array}\right).
\end{gather*}
Here $\EEE_0=\OOO_E\oplus\OOO_E(-\infty)$, the f\/irst elm applied
to $\EEE_0$ gives $\EEE_1=\OOO_E\oplus\OOO_E(p_+ -\infty)$, and
the second one transforms $\EEE_1$ into a stable vector bundle
$\EEE_2$ which f\/its into the exact triple
\[
0\ra\OOO_E \ra\EEE_2\ra\OOO_E(q_1+p_+ -\infty)\ra 0.
\]
Thus, the resulting vector bundle $\EEE=f_{*}(\LLL
)=\EEE_2(-\infty)$ behaves exactly as in the general case
($\tilde{p_+}\neq\tilde{q_2})$.
\end{proof}

Next we will discuss Gabber's elementary transforms as def\/ined
by Esnault and Viehweg~\cite{EV-2}. Gabber's transform of a pair
$(\EEE, \nabla)$, consisting of a vector bundle $\EEE$ over a
curve  and a logarithmic connection on $\EEE$ is another pair
$(\EEE', \nabla')$, where $\EEE'$ is an elementary transform of
$\EEE$ at some pole~$p$ of~$\nabla$, and one of the eigenvalues of
$\Res_p \nabla'$ dif\/fers by~$1$ from the respective eigenvalue
of $\Res_p \nabla$, whilst the other eigenvalues as well as the
other residues remain unchanged. We adapt the def\/inition of
Esnault--Viehweg to the rank-2 case and to our notation:

\begin{definition}\label{def-gabber}
Let $\EEE$ be a rank-$2$ vector bundle on a curve $X$, $\nabla$ a
logarithmic connection on $\EEE$, $p\in X$ a pole of $\nabla$, and
$v\in \EEE|_p$ an eigenvector of the residue
$\Res_p(\nabla)\in\End (\EEE|_p)$. The Gabber transform
$elm_{p,v}(\EEE,\nabla)$ is a pair $(\EEE', \nabla')$ constructed
as follows:
\begin{enumerate}\itemsep=0pt
\item[(i)] $\EEE' =elm_{p,v}^{+}(\EEE)$.

\item[(ii)] $\nabla'$ is identif\/ied with $\nabla$ under the
isomorphism $\EEE|_{X-p}\simeq\EEE'|_{X-p}$ as a meromorphic
connection over $X-p$, and this determines $\nabla'$ as a
meromorphic connection over $X$.
\end{enumerate}
\end{definition}

By a local computation of $\nabla'$ at $p$ one proves:

\begin{lemma} In the setting of Definition {\rm \ref{def-gabber}}, let us
complete $v$ to a basis $(e_1=v,e_2)$ of $\EEE$ near $p$, so that
$\EEE'_p=\OOO_p\cdot\frac{1}{\tau_p}v +\OOO_p\cdot e_2$ and the
matrix $R$ of $\Res_p(\nabla)$ has the form
$R=\left(\begin{array}{cc} \lambda_1& *\\
0&\lambda_2\end{array}\right)$. Then $\nabla'$ is a logarithmic
connection on $\EEE'$ and the matrix $R'$ of its residue at $p$
computed with respect to the basis
$(e'_1,e'_2)=(\frac{v}{\tau_p},e_2)$ of $\EEE'$ has the form
$R'=\left(\begin{array}{cc} \lambda_1-1& 0\\
*&\lambda_2\end{array}\right)$.
\end{lemma}

\begin{theorem}[Bolibruch--Esnault--Viehweg \cite{AB,EV-2}] \label{thm-EV}
Let $\EEE$ be a rank-$r$ vector bundle on a~curve $X$, $\nabla$  a
logarithmic connection on $\EEE$, and assume that the pair $(\EEE,
\nabla)$ is irreducible in the following sense: $\EEE$ has no
$\nabla$-invariant subbundles $\FFF \subset\EEE$. Then there
exists a sequence of Gabber's transforms that replaces $(\EEE,
\nabla)$ by another pair $(\EEE', \nabla')$, in which $\EEE'$ is a
semistable vector bundle of degree $0$ and $\nabla'$ is a
logarithmic connection on $\EEE'$ with the same singular points
and the same monodromy as $\nabla$.
\end{theorem}

We are illustrating this theorem by presenting explicitly one
elementary Gabber's transform which transforms our bundle
$\EEE=f_*\LLL$ of degree $-1$ into a semistable bundle $\EEE'$ of
degree 0:

\begin{proposition} \label{illu}
Let $\EEE$, $\nabla$  be as in Proposition {\rm
\ref{generic-case}}. Let $v$ be an eigenvector of
$\Res_{p_{+}}(\nabla)$ with eigenvalue $\frac{1}{2}$ (see formula
\eqref{Residues}). Then the Gabber transform
$(\EEE',\nabla')=elm^{+}_{p_{+}, v} (\EEE,\nabla)$ satisfies the
conclusion of the Bolibruch--Esnault--Viehweg theorem: $\EEE'$ is
semistable of degree $0$ and $\nabla'$ is a~logarithmic connection
with the same singularities and the same monodromy as $\nabla$.
Furthermore, $\EEE'\simeq\OOO_E (p_+ -\infty)\oplus\OOO_E (q_1
+q_2 -2\infty)$.
\end{proposition}

\begin{proof}
By Corollary \ref{two-elms}, $\EEE'$ is the result of application
of three positive elms to
$\EEE_0(-\infty)=\OOO_E(-\infty)\oplus\OOO_E(-2\infty)$:
\[\EEE'=elm^{+}_{p_{+}, v}
elm^{+}_{q_1, v_2^{(1)}}elm^{+}_{q_2, v_2^{(2)}}
(\EEE_0(-\infty)).
\]
The surface $S_0=\PP(\EEE_0(-\infty))$ can be decomposed as the
open subset $S_0\setminus \bar s_1$ (see Fig.~\ref{Fig.7.}), which
is a line bundle over $E$ with zero section $\bar s_2$, plus the
``inf\/inity section'' $\bar s_1$. The line bundle is easily
identif\/ied as the normal bundle to $\bar s_2$ in $S_0$:
$S_0\setminus \bar s_1\simeq \NNN_{\bar s_2/S_0}\simeq
\OOO_E(\infty)$. Then the pencil $\lvert \bar s_2\rvert= \lvert
\bar s_1+f_{\infty}\rvert$ is the projective line which naturally
decomposes into the af\/f\/ine line $H^0(E,\OOO(\infty))$ and the
inf\/inity point representing the reducible member of the pencil
$\bar s_1+f_\infty$ (the curves $\bar s_{P_1}$, $\bar s_{P_2}$
shown on Fig.~\ref{Fig.7.} are members of this pencil). The fact
that all the global sections $s\in H^{0}(\OOO_E (\infty))$ come
from $H^{0}(\OOO_E)$=\{constants\}\ under the embedding
$\OOO_E\hookrightarrow\OOO_E (\infty)$ implies that they all
vanish at $\infty$. Thus all the $\bar s\in \lvert \bar s_2\rvert$
pass through the point $f_{\infty}\cdot\bar s_2$ which is the zero
of the f\/iber of the line bundle $\OOO_E (\infty)$ over $\infty$.

Using this representation of $S_0$, we can prove the existence of
a cross-section $\bar r\subset S_0$, $\bar r\in \lvert \bar
s_2+f_{p_{+}}\rvert$ passing through the three points $P_0=[p_{+},
v]$ and $P_i=[q_i, v_2^{(i)}]$ ($i=1,2$). Namely, the curves from
the linear system $\lvert \bar s_2+f_{p_{+}}\rvert$ are the
sections of $\OOO_E (\infty+p_+)$ considered as sections of
$\OOO_E(\infty)$ having a simple pole at $p_+$. The fact that they
have a simple pole at $p_+$ means that they meet $\bar s_1$ at
$f_{p_+}\cdot \bar s_1$. The vector space
$H^{0}(\OOO_E(\infty+p_+))$ is $2$-dimensional, so we can f\/ind
$r$ in it taking the values $v_2 ^{(1)}$, $v_2^{(2)}$ at $q_1$,
resp.~$q_2$.

We have $\bar s_1^{2}=-1$, $\bar s_2^{2}=1$, $\bar r^{2}=3$, $\bar
s_1\cdot \bar r=1$, and $\bar s_1\cap \bar r=P_0$. After we
perform the $3$ elementary transforms at $P_i$ ($i=0, 1, 2$), the
self-intersection $\bar r^2$ goes down by $3$. At the same time
$\bar s_1^{2}$ goes up by $2$ when making elms $P_1$, $P_2$ and
descends by $1$ after the elm at $P_0$. Hence in $S'=\PP(\EEE')$,
we have two disjoint sections $\bar r$, $\bar s_1$ with
self-intersection $0$. Thus, by Lemma~\ref{direct-sum},
$\EEE'=\LLL_1\oplus \LLL_2$, where $\LLL_1$, $\LLL_2$ are line
bundles of the same degree. By Lemma~\ref{elm-on-subb},
$\deg\EEE'=\deg\EEE+3=0$, hence $\deg\LLL_1=\deg\LLL_2=0$. The
direct sum of line bundles of the same degree is strictly
semistable.

Next, $\bar s_1$ (in $S_0$) corresponds to the line subbundle
$\OOO_E (-\infty)$. It remains $\OOO_E (-\infty)$ after elms in
$P_1, P_2$, and becomes $\OOO_E(p_+ -\infty)$ after the elm in
$P_0$. Hence $\LLL_1=\OOO_E(p_+ -\infty)$ and
$\LLL_2=\det\EEE'\otimes L^{-1}$. But
$\det\EEE'=\det\EEE(q_1+q_2+p_+)=\OOO_E (q_1+q_2+p_+ -3\infty)$.
Thus $\LLL_2=\OOO_E(q_1+q_2-2\infty)$.
\end{proof}

\begin{remark}\label{0and-1}
 If we f\/ix $E$ and let vary $p_+$, $q_1$, $q_2$, then we see that the generic direct sum $\LLL _1\oplus \LLL _2$ of two line bundles of degree $0$ occurs as the underlying vector bundle of $\nabla'$. According to~\cite{Tu}, the moduli space of semistable rank-$2$ vector bundles on $E$ is isomorphic to the symmetric square~$E^{(2)}$ of $E$, and its open set parameterizes, up to an isomorphism, the direct sums $\LLL _1\oplus \LLL _2$. Thus we obtain a
natural map from the parameter space of our direct image
connections to the symmetric square~$E^{(2)}$, whilst using
the {\em stable} bundles $\EEE$ of degree $-1$ provides a natural
map onto $E$ (Corollary \ref{map-to-E}).
\end{remark}

\begin{remark} Korotkin \cite{Kor-1} considers twisted rank-$2$ connections on $E$ with connection matrices~$A$ satisfying the transformation rule
\begin{gather}\label{eqa} T_a(A)=QAQ^{-1},\qquad T_b(A)=RAR^{-1} \end{gather} for some $2\times 2$ matrices $Q$, $R$.
In the case  when $Q$, $R$ commute, such a twisted connection can be
understood as an ordinary connection on a nontrivial vector bundle
$\EEE$ over $E$ that can be described as follows: let
$E=\CC/\Lambda$ where $\Lambda$ is the period lattice of $E$ with
basis $(1, \tau)$, and let $z$ be the f\/lat coordinate on $E$ (or
on the universal cover $\CC$ of $E$) such that $T_a(z)=z+1$,
$T_b(z)=z+\tau$. Let us make $\Lambda$ act on $\CC^{2}\times\CC$
by the rule
\[
 (v,z)\stackrel{a}\mapsto(Qv, z+1), \qquad (v,z) \stackrel{b}\mapsto (Rv, z+\tau) .
\]
  Then $\EEE\ra E$ is obtained as the quotient $\CC^{2}\times\CC/\Lambda\ra\CC/\Lambda$ of the trivial vector bundle $\CC^{2}\times\CC\xymatrix@1{\ar[r]^{pr_{2}}&}\CC$.

However, the twisted connections obtained in \cite{Kor-1} satisfy
(\ref{eqa}) with non-commuting $Q$, $R$, given by Pauli matrices:
\[ Q=\sigma_1=
\left( \begin{array}{cc} 0 & 1  \\ 1 &  0
\end{array}\right),\qquad R=\sigma_3=\left( \begin{array}{cc} 1 &
0 \\ 0 & -1
\end{array}\right).
\]

This follows from the relation $A=\ud\Psi\Psi^{-1}$, where $\Psi$
is a fundamental matrix of the connection, and from the
transformation law for $\Psi$: $T_a(\Psi)=i\sigma_1\Psi$,
$T_b(\Psi)=i\sigma_3\Psi e^{-2i\pi\lambda\sigma_3}$, where
$\lambda\in\CC$ is a parameter (see (3.74) in loc. cit). Hence
Korotkin's connections are really twisted and have no underlying
vector bundles. This is a major dif\/ference between the result of
\cite{Kor-1} and that of the present paper. Another dif\/ference,
concerning the method, is that the starting point in \cite{Kor-1}
is an ad hoc expression for $\Psi$ in terms of Prym theta
functions of the double cover $C\ra E$, and the connection matrix
is implicit.
\end{remark}

\section[Monodromy and differential Galois groups]{Monodromy and dif\/ferential Galois groups}\label{monodromy}

Let $G$ be the  monodromy group of the connection $\nabla_\EEE $
on $f_{*}\LLL =\EEE$ def\/ined by formula (\ref{Conn_Mat}). It is
the subgroup of $GL(2,\CC)$ generated by $M_a$, $M_b$,
$M_{\gamma_1}$. We will f\/irst consider the case of generic
values of the parameters
$(\lambda_1K,\lambda_1L,\lambda_2K',\lambda_2L')$. Here, {\em
generic} means that the point belongs to the complement of a
countable union of af\/f\/ine $\QQ$-subspaces of $\CC^4$. More
exactly, we require that the triples $(i\lambda_2 K', \lambda_2
L', i\pi)$ and $(\lambda_1 K, i\lambda_1 L, \pi)$ are free over
$\QQ$. Let
\[
{R^{\theta}}=\left( \begin{array}{cc} \cos\theta & i\sin\theta
\\   i\sin\theta &  \cos\theta
\end{array}\right)
, \qquad {H^{\theta}}=\left( \begin{array}{cc} \cosh\theta &
\sinh\theta   \\   \sinh\theta  &  \cosh\theta
\end{array}\right)\qquad (H^{i\theta}=R^{\theta})
.\]

Let $N$ be the normal subgroup of $G$ def\/ined by
\begin{gather}\label{N-via-det} N=\{X\in G\ \mid\  \det X =\pm 1\}.\end{gather}
We have:
\begin{gather}\label{N-general} N=\left\{{\prod_{i=1}^{r} M^{j_i}_{a} M^{k_i}_{b} M^{\epsilon_i}_{\gamma_1}\mid r\geq 0, j_i\in\ZZ, k_i\in\ZZ, \epsilon_i\in \{0,1\}},\sum_{i} k_{i}=\sum_{i}j_{i}=0\right\}.
\end{gather}
We can write  $G$  as the semi-direct product of $N$ with the
subgroup of $G$ generated by~$M_a$,~$M_b$. The latter is
identif\/ied with $\ZZ\times\ZZ$, so $G=N\rtimes(\ZZ\times\ZZ)$.
Let $N_1$ be the subgroup of $N$ generated by $R^{4 {\lambda_1}
K}$, $H^{4 {\lambda_1}L}$. As $M_a=e^{-2i\lambda_2 K'}
R^{-2\lambda_1 K}$, and $M_b=e^{2\lambda_2 L'} H^{2\lambda_1 L}$,
we have $[M_a, M_{\gamma_1}]=R^{-4\lambda_{1} K}$, $[M_b,
M_{\gamma_1}]=H^{4\lambda_{1} L}$, $[M_a, M_b]=1.$  Hence $N$ is
the semi-direct product $N=N_1\rtimes\mu_2$, where $\mu_n\simeq
\ZZ/n\ZZ$ denotes a cyclic group of order $n$, and the factor
$\mu_2$ of the semi-direct product is generated by $M_{\gamma_1}$.
Finally, we obtain a normal sequence $1\lhd N_1\lhd N \lhd G$ with
successive quotients $\ZZ\times\ZZ$, $\ZZ/2\ZZ$, $\ZZ\times\ZZ$,
all of whose levels are semi-direct products. We can write:
\[
G\simeq((\ZZ\times\ZZ)\rtimes\ZZ/2\ZZ)\rtimes(\ZZ\times\ZZ).
\]

We have also $N_1=D(G)$, the commutator subgroup of $G$. As
$D(G)\simeq\ZZ\times\ZZ$ is Abelian, $G$~is solvable of height
$2$.

From now on, we go over to the general case. The formulas
(\ref{N-via-det}), (\ref{N-general}) are no more equivalent. Let
us def\/ine $N$ by (\ref{N-general}), and $N_1$ by the same
formula with the additional condition $\sum_{i}  \epsilon_i\equiv
0 (2).$ We have again the normal sequence $1\lhd N_1\lhd N \lhd
G$. Its f\/irst level is a semidirect product,
$N=N_1\rtimes\mu_2$, but the upper one may be a nonsplit
extension. Def\/ine two group epimorphisms
\begin{alignat}{3}
& \ZZ\times\ZZ  \xymatrix@1{\ar[r]^{\phi_{1}}&}  N_1   ,\qquad  &&
\ZZ\times\ZZ  \xymatrix@1{\ar[r]^{\phi_{2}}&}  G/N ,&\nonumber \\
& (n_1, n_2)  \longmapsto  \frac{\sigma (n_1,n_2)^2}{\det\sigma
(n_1,n_2)},  \qquad &  & (n_1, n_2)  \longmapsto   \sigma
(n_1,n_2)N , & \label{phi_i}
\end{alignat}
where $\sigma (n_1,n_2)=M^{n_1}_{a} M^{n_2}_{b}$.

Thus both $N_1$ and $G/N$ are quotients of $\ZZ\times\ZZ$. We want
to f\/ind out, which pairs $Q_1$, $Q_2$ of quotients of
$\ZZ\times\ZZ$ can be realized as the pair $N_1$, $G/N$ for some
monodromy group $G$. We will denote by $\pi_N$ the canonical
epimorphism $G\rar G/N$, and the maps $\bar\phi_1$, $\bar\phi_2$
are def\/ined by the following commutative diagram:
\[
\xymatrix{ & \ZZ\times\ZZ \ar @{->>} _{\phi_2} [dl]    \ar @{->>}
^{\sigma} [d]
                    \ar  @{->>} ^{\phi_1} [dr]  &      \\
G/N & \langle  M_a,M_b\rangle  \ar @{->>} _{\bar\phi_2} [l]
                            \ar @{^{(}->}  [d]
                               \ar  @{->>} ^{\bar\phi_1} [r]  &    N_1      \\
& G \ar  @{->>} _{\pi_N} [lu]   & \\
}
\]
One can also give $\bar\phi_1$ by the formulas
\[
\bar\phi_1(X)=\frac{1}{\det X}X^2=[X, M_{\gamma_1}] \qquad
\mbox{for all}\ \ X\in \langle  M_a,M_b\rangle .
\]

\begin{proposition}\label{towers}
For any connection \eqref{Conn_Mat}, its monodromy group $G$ fits
into a normal sequence $N_1\lhd N \lhd G$ in such a way, that the
following properties are verified:
\begin{enumerate}\itemsep=0pt
\item Both $N_1$ and $G/N$ are quotients of $\ZZ\times\ZZ$, and
$N/N_1\simeq\mu_2$. \item The extension $N_1\lhd N$ is always
split: $N\simeq N_1\rtimes \mu_2$, the generator $h\in \mu_2$
acting on $N_1$ via the map $g\mapsto g^{-1}$. \item The subgroup
$\langle  M_a,M_b\rangle $ of $G$ provides a splitting of the
extension $N \lhd G$ if and only if $\bar\phi_2$ is an
isomorphism. In this case, the action of $G/N$ on $N$ defining the
split extension is given by $x:g\mapsto g$ and $x:h\mapsto
\bar\phi_1\bar\phi_2^{-1}(x)h$ for any $x\in G/N$, $g\in N_1$.
\end{enumerate}

Conversely, let $(Q_1,Q_2)$ be a pair of group quotients of
$\ZZ\times\ZZ$. Then $(Q_1,Q_2)$ can be realized as the pair
$(N_1,G/N)$ for the monodromy group $G$ of a connection
\eqref{Conn_Mat} if and only if $(Q_1,Q_2)$ occurs in the
following table:
\medskip

\centerline{{\em
\begin{tabular}{|c|c|c|c|c|c|}
\hline N$^\circ$ & $\rk Q_1$ & $\rk Q_2$ & $Q_1$ & $Q_2$ & Restrictions\\
\hline 1$^*$ & 2 & 2 &  $\ZZ\times\ZZ$ &  $\ZZ\times\ZZ$ & ---\\
2 &2 & 1 & $\ZZ\times\ZZ$ & $\mu_{d}\times\ZZ$ & $2|d$ \\
3 &2 & 0 & $\ZZ\times\ZZ$ & $\mu_{2}\times\mu_{d}$ & $2|d$ \\
\hline 4$^*$ & 1 & 2 & $\mu_{d}\times\ZZ$ & $\ZZ\times\ZZ$ & $d\geq 1$ \\
5 & 1 & 1 & $\mu_{d}\times\ZZ$ &$\mu_{d'}\times\ZZ$ & {\em if $2|d$, then} $2|d'$ \\
6 &1 & 0 & $\mu_{d}\times\ZZ$ & $\mu_{d'}$ & $2\nmid d$, $2|d'$ \\
7 &1 & 0  & $\mu_{d}\times\ZZ$ & $\mu_{2}\times\mu_{d'}$ & $d\geq 1$, $2|d'$ \\
\hline 8$^*$ & 0 & 2 & $\mu_{d}$ & $\ZZ\times\ZZ$ & $d\geq 1$ \\
9 & 0 & 1 & $\mu_{d}$ & $\mu_{d'}\times\ZZ$ & $d\equiv d'\!\!\mod\! 2$ \\
10 & 0& 0 &$\mu_{d}$ & $\mu_{d'}$ & $d\geq 1$, $d'\geq 1$ \\
11 & 0& 0&$\mu_{d}$ & $\mu_{2}\times\mu_{d'}$ & $2|d$, $2|d'$ \\
\hline
\end{tabular}
}}\medskip

\noindent The items whose numbers are marked with an asterisk
correspond to the pairs that always give a split extension $N \lhd
G$.
\end{proposition}

\begin{proof}
The f\/irst part, resuming the properties of the tower of group
extensions \mbox{$N_1\lhd N \lhd G$}, is an easy exercise, and we
go over to the second one. Given a pair $(Q_1,Q_2)$, we f\/ind out
whether it is possible to choose epimorphisms
$\ZZ\times\ZZ\xymatrix@1{\ar[r]^{\phi_{1}}&}Q_1$ and
$\ZZ\times\ZZ\xymatrix@1{\ar[r]^{\phi_2} &} Q_2$ and identify them
as the morphisms def\/ined in (\ref{phi_i}) for a suitable choice
of matrices $M_a$, $M_b$. The proof follows a~case by case
enumeration of dif\/ferent types of kernels of $\phi_1$ and
$\phi_2$. To shorten the notation, let us write
$M_a=e^{\alp_1}H^{\beta_1}$, $M_b=e^{\alp_2}H^{\beta_2}$. The case
$\rk\ker\phi_1=\rk\ker\phi_2=0$, corresponding to
$\rk_\QQ(\alp_1,\beta_1,\pi i)=\rk_\QQ(\alp_2,\beta_2,\pi i)=3$,
has been treated before the statement of the proposition. It gives
item 1 of the table.

The proofs of all the other cases resemble each other, and we will
give only one example of this type of argument, say, when both
kernels are of rank~1. Under this assumption, there exist $(d,k_1,
k_2)\in\ZZ^{3}$ and $(d',k'_1, k'_2)\in\ZZ^{3}$ such that
\begin{gather}\label{gcd}
d\geq 1, \qquad d'\geq 1, \qquad \mbox{$\gcd(k_1, k_2)=1$}, \qquad
\mbox{$\gcd(k'_1, k'_2)=1$},
\end{gather}
and
\[
\ker\phi_1=\mbox{$\langle d(k_1, k_2)\rangle $},\qquad
\ker\phi_2=\mbox{$\langle d'(k'_1, k'_2)\rangle $}.
\]
For $(n_1,n_2)\in\ZZ^2$, we have:
\begin{gather}\label{kerphi1}
(n_1,n_2)\in\ker\phi_1 \ \ \equi\ \ \exists\ m\in\ZZ\ | \
n_1\beta_1 +n_2\beta_2=\pi im ;
\\
\label{kerphi2} (n_1,n_2)\in\ker\phi_2 \ \equi\ \exists\
(m_1,m_2)\in\ZZ^2\  \mid  \ M_a^{n_1}M_b^{n_2}=\big(
H^{2\beta_1}\big)^{m_1} \big( H^{2\beta_2}\big)^{m_2}.
\end{gather}
The latter equality can be written in the form $e^\alpha
H^\beta=1$, where
\[
\alpha = n_1\alp_{1}+n_2\alp_2, \qquad \beta =
(n_1-2m_1)\beta_1+(n_2-2m_2)\beta_2.
\]
As $e^\alpha H^\beta=1$ if and only if $e^\alpha =H^\beta=\pm 1$,
we see that the condition of (\ref{kerphi2}) is equivalent to the
existence of an integer vector $(m_0,m_1,m_2,m_3)\in\ZZ^4$ such
that
\begin{gather} m_0\equiv m_3\!\!\!\mod 2
\label{m0m3}\\
n_1\alp_{1}+n_2\alp_2= \pi i m_0, \label{kerphi2-1}\\
(n_1-2m_1)\beta_1+(n_2-2m_2)\beta_2= \pi i m_3 . \label{kerphi2-2}
\end{gather}

Substituting the generators of $\ker\phi_i$ for $(n_1,n_2)$, we
obtain the following system of equations:
\begin{gather}
dk_1\beta_1+dk_2\beta_2=\pi im, \label{m}\\
d'k'_1\alp_{1}+d'k'_2\alp_2= \pi i m_0, \label{m0}\\
(d'k'_1-2m_1)\beta_1+(d'k'_2-2m_2)\beta_2= \pi i m_3 . \label{m3}
\end{gather}
The condition that $d(k_1, k_2)$, $d'(k'_1, k'_2)$ are not just
elements of the corresponding kernels, but their generators, is
transcribed as follows:
\begin{gather}\label{minimality}
\gcd(m,d)=\gcd(d',m_0,2m_1,2m_2,m_3)=1
\end{gather}
for any $(m_1,m_2,m_3)$ satisfying (\ref{m0m3}), (\ref{m3}).

As $\rk\ker\phi_1=1$, the equations (\ref{m}) and (\ref{m3}) have
to be proportional. If $d'$ is odd, but $d$ is even, then at least
one of the coef\/f\/icients of $\beta_i$ in (\ref{m3}) is odd. But
both coef\/f\/icients in (\ref{m}) are even, and this contradicts
(\ref{minimality}). We get the restriction from item 5 of the
table: if $d$ is even, then $d'$ is even, too. This leaves
possible three combinations of parities of $d$, $d'$, and it is easy
to see that a solution to (\ref{gcd}),
(\ref{m})--(\ref{minimality}) exists for any of them. For example,
if $d\equiv d'\!\!\mod\! 2$, then we can choose $k_i$, $k'_i$ in such
a way that $k_i\equiv k'_i\!\!\mod\! 2$ ($i=1,2$), $k_1k'_1\neq
0$. We get a solution to the problem as follows:
\begin{gather*}
m_i=\frac{1}{2}d(k'_i-k_i)\quad (i=1,2),\qquad m=m_0=m_3=1,\qquad
\alp_2=\beta_2=1,\\ \alp_1=\frac{\pi i-d'k'_2}{d'k'_1},\qquad
\beta_1=\frac{\pi i-dk_2}{dk_1}.
\end{gather*}
Our choice for $\alp_2$, $\beta_2$ is explained by the observation
that we should have $\rk_\QQ(\alp_1,\beta_1,\pi
i)=\rk_\QQ(\alp_2,\beta_2,\pi i)=2$, and 1 is the simplest complex
number which is not a rational multiple of~$\pi i$.
\end{proof}

\begin{remark} In the above proof,
if $\ker\phi_2\not\subset\ker\phi_1$, then any solution of
(\ref{m})--(\ref{m3}) satisf\/ies the condition $(m_1,m_2)\neq 0$,
which means that $d'(k'_1, k'_2)\not\in\ker\sigma$. Hence
$\sigma(d'(k'_1, k'_2))$ is a nonzero element of $\ker\bar\phi_2$,
and $\bar\phi_2$ is not an isomorphism. This implies that the
extension $N\lhd G$ is nonsplit. Hence it is never split, unless
$d|d'$. In this case, it can be occasionally split, if
$\ker\phi_2\subset\ker\phi_1$.
\end{remark}

We can deduce from Proposition \ref{towers} a description of all
the f\/inite monodromy groups; they correspond to lines 10 and 11
of the table. This description is only partial, because we do not
determine completely the extension data.

\begin{corollary}
All the finite monodromy groups $G$ of connections
\eqref{Conn_Mat} are obtained as extensions
\[
D_d\hookrightarrow G\twoheadrightarrow \mu_{d'}\qquad  (d\geq 1,\
d'\geq 1)
\]
or
\[
D_d\hookrightarrow G\twoheadrightarrow \mu_2\times\mu_{d'}\qquad
(2\mid d,\ 2\mid d'),
\]
where $D_d=\mu_d\rtimes\mu_2$ is the dihedral group.
\end{corollary}

\begin{corollary}
The only finite Abelian groups occurring as the monodromy groups
of connections~\eqref{Conn_Mat} are $\mu_2$ and
$\mu_2\times\mu_d$\ ($d\geq 2$).
\end{corollary}

We add a few examples of inf\/inite monodromy groups with
nongeneric parameters $(\lambda_1K,\lambda_1L$,
$\lambda_2K',\lambda_2L')$.

\begin{example}
It is easy to select the parameters to get for $G$ one of the
groups $D_n\times \ZZ^i$ or $D_n\rtimes \ZZ^i$, where
$n\in\NN\cup\{\infty\}$, $i=0,1,2$. For example, to get
$D_n\rtimes\ZZ $, we can set $M_a=R^{\frac{2\pi}{n}}$, $M_b=R^1$,
and to get $D_n\times\ZZ $, we can set $M_a=R^{\frac{2\pi}{n}}$,
$M_b=1$.
\end{example}

Now that we have described the structure of the monodromy group of
$\nabla_\EEE $, we can ask the question on its Zariski closure.
According to \cite[Proposition 5.2]{Ka}, the Zariski closure of
$G$ is the dif\/ferential Galois group $\DGal(\nabla_\EEE )$. For
the reader's convenience, we recall its def\/inition.

Let $(K,\; ')$ be a dif\/ferential f\/ield with f\/ield of
constants $\CC$. This means that $K$ is endowed with a
$\CC$-linear derivation \ $':K\rar K$.

\begin{definition}
Let $(K,\;  ')\subset (L,\; ')$ be an extension of dif\/ferential
f\/ields with f\/ield of constants $\CC$. The dif\/ferential
Galois group $\DGal(L/K)$ is the group consisting of all the
$K$-automorphisms $\sigma$ of $L$ such that
$\sigma(f')=(\sigma(f))'$ for all $f\in L$.
\end{definition}

If $L$ is f\/initely generated as a $K$-algebra, say, by $p$
elements, then $\DGal(L/K)$ can be embedded onto $GL(p,\CC)$, and
it is an algebraic group if considered as a subgroup of
$GL(p,\CC)$ in this embedding.

We apply this def\/inition to $K=\CC(E)$, the derivation $'$ being
the dif\/ferentiation with respect to some nonconstant function
$z\in K$. Given a connection $\nabla_\EEE $ on $E$, we can
consider a fundamental matrix $\Phi$ of its solutions, and set $L$
to be the f\/ield generated by all the matrix elements of $\Phi$.
The group $\DGal(\nabla_\EEE )$ is def\/ined to be $\DGal(L/K)$.
See \cite{VDP,VS} for more details.

Remark that the monodromy group $G$ lies in the subgroup $\GG$ of
$GL(2, \CC)$ def\/ined by
\begin{displaymath}
\GG=\left\{\left( \begin{array}{cc} C\alpha & C\epsilon\beta   \\
C\beta  &  C\alpha\epsilon
\end{array}\right)\mid C\in\CC^{*}, (\alpha,\beta)\in\CC^2 , \epsilon\in\{-1,1\}, \alpha^{2}-\beta^{2}=1  \right\} .
\end{displaymath}
Denote by $\GG_0$ the connected component of unity in $\GG$,
singled out by the condition $\epsilon=1$. The Zariski closure
$\bar G$ of $G$ is contained in $\GG$ and is not contained in
$\GG_0$. The following statement is obvious:

\begin{lemma}
 Let $\psi:\CC^{*}\times\CC^{*}\rtimes\{-1,1\}\lra\GG$ be defined by
\[
(\lambda,\mu,\epsilon)\longmapsto \left( \begin{array}{cc}
\lambda\alpha & \lambda\beta \epsilon  \\ \lambda\beta  &
\lambda\alpha\epsilon
\end{array}\right)
\]
with $\alpha=\frac{1}{2} (\mu+\frac {1}{\mu})$, $\beta=\frac{1}{2}
(\mu-\frac {1}{\mu})$. Then $\psi$ is a surjective morphism with
kernel $\{(1, 1, 1)$, $(-1, -1, -1)\}$.
\end{lemma}

We see that $\GG_0=\psi(\CC^{*}\times\CC^{*})$ is identif\/ied
with the quotient $\CC^{*}\times\CC^{*}/\{-1,1\}$, and the latter
is isomorphic to $\CC^{*}\times\CC^{*}$ via the map $(z_1,
z_2)\!\!\!\mod\!\!\{-1,1\}\longmapsto (z_1z_2, \frac{z_1}{z_2})$.
Thus we get an explicit isomorphism $\GG_0\simeq
\CC^{*}\times\CC^{*}$. Using this identif\/ication, one can easily
determine the Zariski closure~$\overline{G_0}$ of the subgroup
$G_0=G\cap \GG_0=\langle M_a,M_b\rangle $ of
$\CC^{*}\times\CC^{*}$, and $\DGal(\nabla_\EEE
)=\overline{G_0}\rtimes \langle  M_{\gamma_1}\rangle $.

We can use the following observations:
\begin{enumerate}\itemsep=0pt
\item[a)] If a pair $(s,t)\in\CC^{*}\times\CC^{*}$ is such that
$\rk_\QQ (\ln(s), \ln(t), i\pi)=1$ (that is, $s$ and $t$ are roots
of unity), then the group generated by the pair $(s,t)$ is
f\/inite and coincides with its closure.

\item[b)] If a pair $(s,t)\in\CC^{*}\times\CC^{*}$ is such that
$\rk_\QQ (\ln(s), \ln(t), i\pi)=2$, and $k_1\ln(s)+k_2\ln (t)+
2k_3 i\pi=0$ is a $\ZZ$-linear relation with relatively prime
$k_i$, then $\overline{\langle (s,t)\rangle }$ is the subgroup $V$
of $\CC^{*}\times\CC^{*}$ def\/ined by $z_1^{k_1}z_2^{k_2}=1$,
isomorphic to $\CC^*\times\mu_d$, where $d=\gcd(k_1,k_2)$, and
$\mu_d$ is the cyclic group of order~$d$.

\item[c)] If the triple $(\ln(s), \ln(t), \pi i)$ is free over
$\QQ$, then the closure of $\langle (s,t)\rangle$ is $\CC^*\times
\CC^{*}$.
\end{enumerate}

Apply this to pairs $(s,t)$ belonging to the subgroup generated by
two pairs $(s_1,t_1)$, $(s_2,t_2)$ which are the images of $M_a$,
resp. $M_b$. Then if $(\ln(s_j), \ln(t_j), \pi i)$ is free over
$\QQ$ for at least one value of $j=1$ or 2, then
$\overline{G_0}=\CC^{*}\times\CC^{*}$ and $\DGal(\nabla_\EEE
)=\GG$. In the case when both triples $(\ln(s_1), \ln(t_1), \pi
i),(\ln(s_2), \ln(t_2), \pi i)$  are not free over $\QQ$, the
necessary and suf\/f\/icient condition for $\overline{\langle
(s_1,t_1),(s_1,s_2)\rangle }$ to be $\CC^{*}\times\CC^{*}$ is the
following: $\rk_\QQ(\ln(s_j), \ln(t_j), i\pi)=2$ for both values
$j=1,2$, and if $a_{j1}\ln(s_j)+a_{j2}\ln(t_j)+ a_{j3}\pi i=0$
($j=1,2$) are nontrivial $\QQ$-linear relations in these triples,
then $\left|\begin{array}{cc}a_{11} & a_{12}\\ a_{21} &
a_{22}\end{array}\right|\neq 0$. This condition can be easily
formulated in terms of the epimorphisms $\phi_i$ def\/ined in
(\ref{phi_i}): $\ker\phi_1$, $\ker\phi_2$ are both of rank 1 and
$\ker\phi_1\cap \ker\phi_2=0$. In this case we have the same
conclusion: $\DGal(\nabla_\EEE )=\GG$.

We obtain the following description of possible dif\/ferential
Galois groups of connections (\ref{Conn_Mat}):

\begin{proposition}
Let $r_i=\rk_\QQ\ker\phi_i$\  $(i=1,2)$.
\begin{enumerate}\itemsep=0pt
\item[(i)] $\DGal(\nabla_\EEE )=\GG$ if and only if one of the
following condition is verified: either $\min \{r_1,r_2\}=0$, or
$r_1=r_2=1$ and $\ker\phi_1\cap \ker\phi_2=0$.

\item[(ii)] $\DGal(\nabla_\EEE )$ is a $1$-dimensional subgroup of
$\GG$ if and only if $\min \{r_1,r_2\}=1$ and the condition of (i)
is not satisfied. Then there exists a one-parameter subgroup $V_0$
and a finite cyclic subgroup $\mu_d$ in $\GG$ such that
$\DGal(\nabla_\EEE )=(V_0\mu_d)\rtimes \langle
M_{\gamma_1}\rangle $.

\item[(iii)] $\DGal(\nabla_\EEE )$ is finite if and only if
$r_1=r_2=2$, and then $\DGal(\nabla_\EEE )=G$.
\end{enumerate}
\end{proposition}

\subsection*{Acknowledgements}
I am greatly indebted to my research advisor D.~Markushevich for
his encouragement and help. I~would like to thank D.~Korotkin for
explaining me some points. I also acknowledge with pleasure the
hospitality of the Mathematics Institute of the Chinese Academy of
Sciences, where was done a part of the work on the article.
The work was partially supported by the Conseil Departemental du Nord.

\pdfbookmark[1]{References}{ref}
\LastPageEnding

\end{document}